\documentclass[oneside,reqno,a4paper,11pt]{amsart}
\usepackage{amsmath,amsfonts,amssymb,}
\usepackage{CJK,bbm,color,soul,supertabular,longtable,verbatim,everb,extarrows}
\usepackage[warning, easyold,math]{easyeqn}
\usepackage{hypernat,hyperref}
\usepackage{graphicx}

\newtheorem {theorem}{Theorem}[section]
\newtheorem {lemma}[theorem]{Lemma}
\newtheorem {corollary}[theorem]{Corollary}
\newtheorem {proposition}[theorem]{Proposition}
\theoremstyle{remark}
\newtheorem {remark}{Remark}[section]

\theoremstyle{problem}
\newtheorem {problem}{Problem}
\theoremstyle{definition}
\newtheorem {definition}{Definition}[section]
\theoremstyle{plain}
\numberwithin {equation}{section}

\makeindex

\begin{document}
\vspace{1cm}

\title[Subsonic Flows in a Multi-Dimensional Nozzle]{Subsonic Flows in a Multi-Dimensional Nozzle$^*$}
\author[Lili Du, \ \ Zhouping Xin,\ \  Wei Yan ]{Lili Du$^{\lowercase{a,b,1}}$, \ Zhouping Xin$^{\lowercase{b,2}}$, \ Wei Yan$^{\lowercase{c,b,3}}$}
\thanks{$^1$ E-Mail: lldumath@hotmail.com. $^2$ E-Mail: zpxin@ims.cuhk.edu.hk. $^3$ E-Mail: wyanmath@gmail.com.}
\maketitle
\begin{center} $^a$ Department of Mathematics, Sichuan Univeristy,

          Chengdu 610064, P. R. China.

           $^b$ The Institute of Mathematical Sciences, The Chinese University of Hong Kong,

           Shatin, NT, Hong Kong

           $^c$ Science and Technology Computation Physics Laboratory,

           Institute of Applied Physics and Computational Mathematics,

           Beijing 100088, P.R. China.
\end{center}

\maketitle

\begin{abstract}
In this paper, we study the global subsonic irrotational flows in a
multi-dimensional ($n\geq 2$) infinitely long nozzle with variable cross sections. The flow is
described by the inviscid potential equation, which is a second
order quasilinear elliptic equation when the flow is subsonic.
First, we prove the existence of the global uniformly subsonic flow
in a general infinitely long nozzle for arbitrary dimension for sufficiently small incoming
mass flux and obtain the
uniqueness of the global uniformly subsonic flow. Furthermore,
we show that there exists a critical value of the incoming mass flux
such that a global uniformly subsonic flow exists uniquely, provided that the incoming mass
flux is less than the critical value. This gives a positive answer
to the problem of Bers on global subsonic irrotational flows in infinitely long nozzles for arbitrary
dimension \cite{BERS}. Finally, under suitable asymptotic assumptions of the nozzle, we obtain the asymptotic behavior of the subsonic flow in
far fields by a blow-up argument. The main ingredients of our analysis are
methods of calculus of variations, the Moser iteration techniques
for the potential equation and a blow-up argument for infinitely long
nozzles.
\end{abstract}
\everymath{\displaystyle}
\newcommand {\eqdef }{\ensuremath {\stackrel {\mathrm {def}}{=}}}

% math symbols
\def\Xint #1{\mathchoice
{\XXint \displaystyle \textstyle {#1}} %
{\XXint \textstyle \scriptstyle {#1}} %
{\XXint \scriptstyle \scriptscriptstyle {#1}} %
{\XXint \scriptscriptstyle \scriptscriptstyle {#1}} %
\!\int}
\def\XXint #1#2#3{{\setbox 0=\hbox {$#1{#2#3}{\int }$}
\vcenter {\hbox {$#2#3$}}\kern -.5\wd 0}}
\def\ddashint {\Xint =}
\def\dashint {\Xint -}
\def\clockint {\Xint \circlearrowright } % GOOD !
\def\counterint {\Xint \rotcirclearrowleft } % Good for Computer Modern !
\def\rotcirclearrowleft {\mathpalette {\RotLSymbol { -30}}\circlearrowleft }
\def\RotLSymbol #1#2#3{\rotatebox [ origin =c ]{#1}{$#2#3$}}

\def\aint{\dashint}

\def\arraystretch{2}
\def\eps{\varepsilon}

\def\s#1{\mathbb{#1}} % set
\def\t#1{\tilde{#1}} %new variables
\def\b#1{\overline{#1}}
\def\N{\mathcal{N}} %Nozzle
\def\M{\mathcal{M}} %Mach number
\def\R{{\mathbb{R}}}

\def\ba{\begin{aligned}}
\def\ea{\end{aligned}}
\def\be{\begin{equation}}
\def\ee{\end{equation}}

\def\bes{\begin{mysubequations}}
\def\ees{\end{mysubequations}}

\def\cz#1{\|#1\|_{C^{0,\alpha}}}
\def\ca#1{\|#1\|_{C^{1,\alpha}}}
\def\cb#1{\|#1\|_{C^{2,\alpha}}}

\def\lb#1{\|#1\|_{L^2}}
\def\ha#1{\|#1\|_{H^1}}
\def\hb#1{\|#1\|_{H^2}}

\def\cin{\subset\subset}
\def\Ld{\Lambda}
\def\ld{\lambda}
\def\ol{{\Omega_L}}
\def\sla{{S_L^-}}
\def\slb{{S_L^+}}
\def\C{\mathbf{C}} %%unit cylinder
\def\cl#1{\overline{#1}}
\def\ra{\rightarrow}
\def\xra{\xrightarrow}
\def\g{\nabla}
\def\a{\alpha}
\def\b{\beta}
\def\d{\delta}
\def\th{\theta}
\def\fai{\varphi}
\def\O{\Omega}
\def\f{\frac}
\def\p{\partial}
\def\disp{\displaystyle}

\def\H{\Theta} %%truncate density function

\section{Introduction}
\label{intro} This paper is devoted to the existence and the
uniqueness of global subsonic flows for the Euler equations for
steady irrotational compressible fluids. Our focus is on the global
subsonic flows in a general multi-dimensional infinite nozzle, which
is an important subject in gas dynamics (see \cite{Bers54a}
\cite{BERS} \cite{Courant}\cite{Fei93}\cite{Lan59}).

Consider the steady isentropic compressible Euler equations
\be\label{euler} \left\{\ba
&\ \text{div} (\rho u) = 0,&\quad\text{ in }\O,\\
&\ \text{div}(\rho u\otimes u) + \g p=0,&\quad\text{ in }\O,
\ea\right.\ee where $\rho$, $u=(u_1, \ldots, u_n)$, $p$ represent
the density, velocity, and the pressure of the fluid, respectively.
Moreover, the pressure $p= p(\rho)$ is a smooth function of $\rho$
and $p'(\rho) > 0$, $p''(\rho)>0$ for $\rho > 0$.

It is easy to derive the following so-call Bernoulli's law
\cite{Courant}
\begin{equation}\label{bernoulli}
u\cdot\g\left(\frac12|u|^2+h(\rho)\right)=0,
\end{equation}
where $h(\rho)$ is the enthalpy, defined by $h(\rho) =
\int_1^\rho\frac{p'(s)}sds$. The relation (\ref{bernoulli}) implies
that the quantity $B(\rho, |u|^2)=\f12|u|^2+h(\rho)$, named
Bernoulli's function, remains constant along the stream line in a
steady isentropic flow.

If, in addition, the flow is assumed to be irrotational, ie. the
vorticity of the flow velocity
$$
\g \times u = 0,\ \ \ \ \ \text{in} \ \ \O,
$$
then there exists a velocity potential function $\fai$, at least
locally, such that
$$
u(x) = \g\fai(x).
$$
In this case, the relation (\ref{bernoulli}) simplifies to the
following strong version of the Bernoulli's law \be\label{strongber}
\g B(\rho, |\fai|^2) = \g
\left(\frac{1}{2}|\g\fai|^2+h(\rho)\right)=0. \ee This yields a
density-speed relation for steady irrotational flows. Therefore, the
density $\rho$ can be determined by the speed $|\g\fai|$, denoted by
$\displaystyle{\rho\left(|\g\fai|^2\right)}$. Then the steady Euler
equations (\ref{euler}) are reduced to the following well-known
scaler potential equation \be\label{potential}
\text{div}\left(\rho(|\g\fai|^2)\g\fai\right)=0,\ \ \ \ \text{in}\ \
\O. \ee

One of the most important parameters to the fluid dynamics is
the {\it Mach number}, which is defined as a non-dimensional ratio of
the fluid velocity to local sound speed,
$$
M=\frac{|u|}{c(\rho)},
$$
where $c(\rho)= \sqrt{p'(\rho)}$ is the local sound speed. Mathematically,
the second-order nonlinear equation (\ref{potential}) is elliptic in the
subsonic region, ie. $M < 1$ and hyperbolic in the supersonic region
where $M > 1$.

Subsonic flows are those in which the local velocity speed is
smaller than sonic speed everywhere, i.e. the Mach number of the
flow is less than $1$. Since the corresponding equations of subsonic
flows possess some elliptic properties, problems related to subsonic
flows are, in general, have extra-smoothness to those related to
transonic flows or supersonic flows. There are many literatures in
this field in the past decades. The first result is due to Frankl
and Keldysh \cite{Frankl34}. They studied the subsonic flows around
a 2D finite body (or airfoil) and proved the existence and the
uniqueness for small data by the method of successive
approximations. Later on, Bers \cite{Bers49}\cite{Bers51} proved the
existence of subsonic flows with arbitrarily high local subsonic
speed for the Chaplygin gas (minimal surface). By a variational
method, Shiffman \cite{Shiffman52b}\cite{Shiffman52} proved that, if
the infinite free stream flow speed $u_\infty$ is less than some
critical speed, there exists a unique subsonic potential flow around
a given profile with finite energy. Shortly afterwards, Bers
\cite{Bers54} improved the uniqueness results of Shiffman. Finn and
Gilbarg \cite{Finn572d} proved the uniqueness of the 2D potential
subsonic flow about a bounded obstacle with given circulation and
velocity at infinity.
All above the results are related to two dimensional problems.
For three (or higher) dimensional case, Finn and Gilbarg
\cite{Finn573d} proved the existence, uniqueness and the asymptotic
behavior with implicit restriction on Mach number $M$. Payne and
Weinberger \cite{Payne57} improved their results soon after. Later,
Dong \cite{Dongguagnchang91} extended the results of Finn and
Gilbarg \cite{Finn573d} to any Mach number $M<1$ and to arbitrary
dimensions. Furthermore, in \cite{Dongguangchang93}, Dong and Ou
extended the results of Shiffman to higher dimensions by the direct
method of calculus of variations and the standard Hilbert space
method.

\begin{figure}[!h]
\includegraphics[width=110mm]{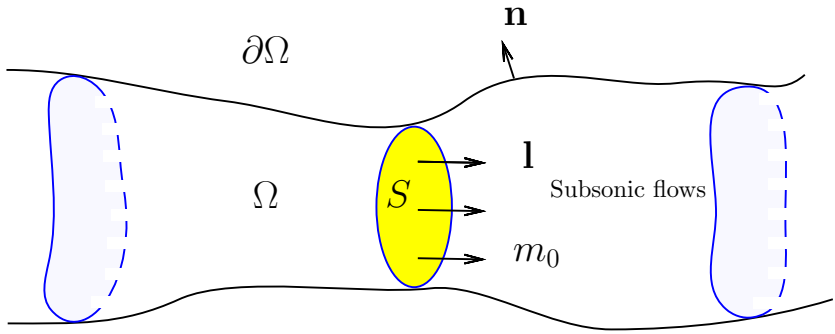}\\
\caption{Subsonic flow in a nozzle} \label{fig1}
\end{figure}

All results as above (including \cite{Shiffman54}-\cite{Gilbarg55})
are related to the subsonic flows past a profile. Another important problem is the study of subsonic flows is the theory of global subsonic flow in a variable nozzles as formulated by Bers in \cite{BERS}:
\begin{problem}\label{op} Find $\fai$ such that,
\be\label{nozzle} \left\{\begin{array}{ll}
 \text{div}\bigg(\rho(|\g\fai|^2)\g\fai\bigg) = 0,&\quad\text{in}\ \ \O,\\
\ \ \f{\p\fai}{\p \vec{n}} = 0,&\quad\text{on}\ \ \p\O,\\
\int_{S_0}\rho\left(|\g\fai|^2\right)\f{\p \fai}{\p \vec{l}}dS = m_0>0,& \\
|\g\fai| < c(\rho),&\quad\text{in}\ \ \O, \end{array}\right. \ee
where $\O\subset \R^n$ is an infinitely long nozzle, $m_0>0$ is the
mass flux passing through the nozzle, $S_0$ is an arbitrary cross
section of the nozzle, $\vec n$ and $\vec l$ are the unit outer
normal of the domain $\O$ and $S_0$, respectively (Please see Fig.
\ref{fig1}).
\end{problem}

In the famous survey \cite{BERS}, Bers claimed without proof
the unique solvability of sufficiently slow subsonic irrotational flows
in two dimensional channel.
The rigorous proof of this fact was achieved mathematically recently by Xie and Xin \cite{XX1}. They established a very
complete, satisfactory and systematic theory for the two dimensional
subsonic flows in an infinitely long nozzle for potential flows, which not only solves the Problem 1 in this case, but also yields the existence of subsonic-sonic flows in the nozzle as limits of subsonic flows. One of the key ideas in \cite{XX1}, is to use the
stream function to formulate the problem to a quasilinear elliptic
problem with Dirichlet boundary conditions. The benefit of the
stream function formation of the problem is that, the stream
function $\psi$ has a priori $L^\infty$ bound, and the flow region
of two dimensional nozzle, though infinitely long, has finite
"width". So one can obtain the boundary $L^\infty$ estimate of the
gradient of the stream function, $\g\psi$, by constructing proper
barrier functions and the standard comparison principle for
subsolution to second order elliptic equation. Similar approach has
been applied in 3D axis-asymmetric nozzles by Xie and Xin in
\cite{XN}. Furthermore, these ideas are also useful to study the physically more important case, subsonic Euler flows, by Xie and Xin in \cite{xx2} (see also the generalization in \cite{dd}). However, it seems difficult to
apply the method in \cite{XX1} and \cite{XN} in general
multi-dimensional ($n\geq 3$) nozzles, since the stream function
formulation can not work in this case. Thus, we have to consider a
different approach from that in \cite{XX1} to treat the subsonic
problem in multi-dimension case.

On the other hand, since the domain of an infinitely long nozzle is differentiable homeomorphism to an infinitely long
cylinder which is unbounded, the nozzle flow problems
are different to the airfoil problems in which the
domains are exterior domains. The main advantage of the exterior
domain is that it can be transformed to a bounded domain through a
Kelvin-like transformation. Then the airfoil problem can be
transformed (explicitly or implicitly) to a scalar quasilinear
elliptic problem with a bounded domain. This feature of the exterior
domain plays an essential role in the previous airfoil
results. For instance, in \cite{OU94}, a Hardy-type inequality in
the exterior domain is essential. But there is no similar Hardy-type
inequality for the domain of nozzle flows, which is the one of the
main difficulties in our case. For more detailed discussions, we refer to \cite{yan}.

The main purpose of this paper is to study subsonic flows in general
multi-dimensional ($n\geq 2$) infinitely long nozzles. First, we
formulate a subsonic truncated problem, which is a uniformly
elliptic equation in a bounded domain. Moreover, we prove the
existence of the weak solution to the truncated problem by a
variational method, and use the approximated variational problems in
bounded domains to approximate the original Problem \ref{op}. To realize this procedure, some uniform
estimates are needed to show that the approximated solutions converge to the ones of
the original Problem \ref{op}. However, one can not expect to get
the uniform boundary gradient estimate of $\fai$ by the classical
barrier function argument, since the potential function $\fai$ is
essentially unbounded, which is another main difficulty in this
paper. The key observation here is that, though the potential
function $\fai$ is unbounded, the $L^2$ average of $\g\fai$ is
uniformly bounded (see the estimate (\ref{L})). Using this
fact and the uniform ellipticity, we prove the "local average
estimate" which states that the average estimate implies the local
average of the gradient $\g\fai$ is uniformly bounded (see
(\ref{eq-LAE}) for details). That is, $\g\fai$ is locally $L^2$
bounded. Then, it is easy to get the $L^\infty$ bound of $\g\fai$ by
the standard Moser iteration. With this key estimate of uniformly
$L^\infty$ bound of $\g\fai$, we establish the existence of the
subsonic flows in an infinitely long nozzle for arbitrary dimensions
for suitable small incoming mass flux, including the two dimensional
case in \cite{XX1}. Next, we show that the global uniformly
subsonic flow is unique. The proof is based on considering the
linear equation satisfied be the difference of two solutions of the
nonlinear potential equation. Moreover, we prove the existence of
the critical incoming mass flux for subsonic flows. Finally, with
the additional asymptotic assumptions on the nozzle at the far
field, we obtain some asymptotic behaviors of the subsonic flow at
the far field by a blow-up argument.

Before stating the main results in this paper, we first give the
following assumptions on the nozzle.

{\bf Basic assumptions on $\O$}. There exists an invertible
$C^{2,\alpha}$ map $T:\ \cl\O\ra\cl\C\ : x\mapsto y$ satisfying
\be\label{H1} \left\{
\begin{array}{l}
T(\p\O) = \p\C,\\
\text{ For any }k\in\R,\ T(\O\cap\{x_n=k\}) = B(0,1)\times\{y_n=k\},\\
\disp \cb T,\ \cb{T^{-1}} \leq K,
\end{array}
\right. \ee where $K$ is a uniform constant, $\C = B(0,1)\times
(-\infty, \infty)$ is a unit cylinder in $\R^n$, $B(0,1)$ is unit
ball in $\R^{n-1}$ centered at the origin, $x_n$ is the longitudinal
coordinate.

\begin{figure}[!h]\centering
\includegraphics[width=90mm]{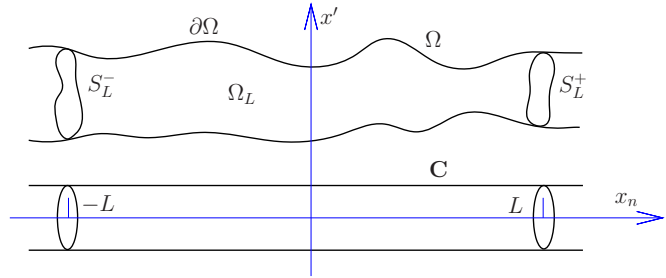}\\
\caption{Basic assumptions on $\O$} \label{fig2}
\end{figure}

{\bf Asymptotic assumptions on $\O$}. Suppose that the nozzle
approaches to a cylinder in the far fields, ie. \be\label{AA} \Omega
\cap \{x_n=k\} \rightarrow S_\pm , \ \ \ \ \ \ \text{as}\ \ \ k
\rightarrow \pm \infty, \ee respectively, where $S_\pm$ are $n-1$
dimensional, simply connected, $C^{2,\alpha}$ domains.

\begin{figure}[!h]\centering
\includegraphics[width=90mm]{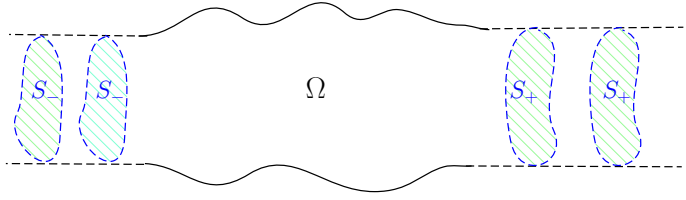}\\
\caption{Asymptotic assumptions on $\O$} \label{fig3}
\end{figure}

{\bf Nondimensionalization of the quantities.} It follows from
Bernoulli's Law (\ref{bernoulli}) that in a potential flow the
density is a given function of speed. Applying the fact
(\cite{BERS}, \cite{Courant}) that there exists a critical speed
$q_{cr}$ such that the flow is subsonic if the speed is less than
$q_{cr}$, we can introduce the nondimensional velocity and density as
$$\hat u=\f u{q_{cr}},\ \ \hat v = \f v{q_{cr}},\ \ \hat \rho=\f\rho{\rho(q_{cr}^2)}.$$
With an abuse of the notation, we still denote the nondimesional
quantities by $u,v,\rho$. Then it is easy to check that $\rho q\leq 1$
for $q\geq 0$ and that the flow is subsonic provided that $q<1$ or
$\rho > 1$.

Our main results in this paper are stated as follows.
\begin{theorem}\label{mainth} Suppose that the nozzle $\O$ satisfies the basic assumptions
(\ref{H1}). Then

(i) there exists a positive number $M_0$ depending only on $\O$,
such that if $m_0\leq M_0$, then there exists a uniformly subsonic
flow through the nozzle, ie., the Problem \ref{op} has a smooth
solution $\fai\in C^\infty(\Omega)$. Moreover,
\[
\|\g\fai(x)\|_{C^{1,\a}(\O)}\leq C m_0,
\]
where $C>0$ is a uniform constant independent of $M_0, m_0$, and
$\fai$.

(ii) There exists a critical mass flux $M_c\leq 1$, which depends only
on $\Omega$, such that if $0\leq m_0 < M_c$, then there exists a unique
uniformly subsonic flow through the nozzle with the following
properties

\be\label{1} Q(m_0)=\sup_{x\in\bar\Omega}|\nabla\varphi|<1, \ee and
$Q(m_0)$ ranges over $[0,1)$ as $m_0$ varies in $[0,M_c)$.

(iii) Furthermore, assume that the nozzle satisfies the asymptotic
assumption (\ref{AA}), then the flow approaches the uniform flows at
the far fields, ie. \be\label{2} \nabla \varphi = (0, \cdots,
q_\pm),\ \ \ \ \ \ \ \text{as}\ \ \ x_n\rightarrow \pm\infty, \ee
respectively, with $q_\pm$ being constants determined uniquely by
$$
\rho(q_\pm^2)q_\pm =\frac{m_0}{|S_\pm|},
$$
here $|S_\pm|$ represents the measure of the domain $S_\pm$,
respectively.
\end{theorem}

\begin{figure}[!h]\centering
\includegraphics[width=90mm]{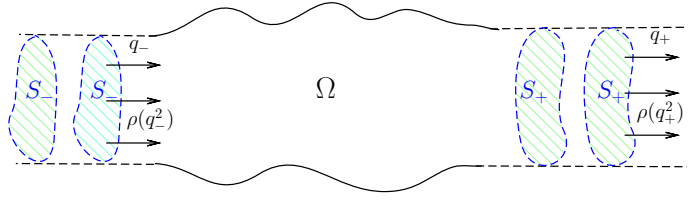}\\
\caption{Asymptotic behaviors of the subsonic flows at the far
fields} \label{fig4}
\end{figure}

\begin{remark}
In the first statement of Theorem \ref{mainth}, it follows from the proof in Section 3 that one can derive an
explicit form of $M_0$, which depends only on the nozzle $\Omega$. In
particular, it does not depend on the equation of the states. On the
other hand, in the second statement of the Theorem \ref{mainth}, we
just give the existence of the critical mass flux $M_c$ for a given
infinite long nozzle. Clearly, $M_0$ is a lower bound of $M_c$.
\end{remark}

\begin{remark} In the proof of the uniqueness of uniformly subsonic flows, it is not necessary to require the asymptotic assumption (\ref{AA}) on the nozzle. It is
quite different from the strategy in \cite{XX1} for the 2D case, in
which the proof of the uniqueness depends on the asymptotic
behaviors of the uniformly subsonic flow in the far fields. However,
in this paper, the uniqueness of the uniformly subsonic flow is
obtained in arbitrary dimensional nozzle without the asymptotic
assumption (\ref{AA}).

\end{remark}

This paper is organized as follows. In the next section, we
introduce some necessary preliminaries. In
Section 3, we prove the first statement of the Theorem 1.1. Our
strategy for the existence of subsonic flows with small incoming mass
flux can be divided into six steps: Step 1, truncate the
coefficients of the potential equation to guarantee the strong
ellipticity and truncate the unbounded nozzle to a series of bounded
domains $\O_L$, to formulate the approximated strong elliptic
problems in bounded domains. Step 2, solve the approximate truncated
problems by a direct variational method. Step 3, improve the
regularity of the variational solutions to give the $H^2$
regularity. Step 4, prove the $L^2$ local average estimates to the
gradient of the solutions. Step 5, obtain the classical $C^{1,\a}$
estimate of the approximate solutions. Step 6, based on these key
estimates, the existence of the subsonic solution to the nozzle
problem for suitable small incoming mass flux is proved. The
uniqueness of the uniformly subsonic solution is given in Section 4,
while the existence of the critical value for incoming mass flux is
obtained in Section 5. In the last section, we prove that the
subsonic nozzle flows approach to the uniform flows at the far
fields when the nozzle satisfies the asymptotic assumption
(\ref{AA}).

In this paper, $x$, $y$ always denote the variables in $\O$ and $\C$
respectively, $\fai$ denotes the function defined in $\O$ and
$\t\fai = \fai\circ T^{-1}$ denotes the corresponding function
defined in $\C$. $\p$, and $\g$ denote the derivatives with respect
to $x$ in $\O$, while $\t\p$, and $\t\g$ denote the derivatives with
respect to $y$ in $\C$.
 $A\sim B$ means
\[
\f 1C A\leq B\leq C A,
\]
with $C$ a positive constant.
\section{Preliminaries}
In this section, we give some basic notations, definitions and facts to be
used in this paper.
\subsection{Morrey theorem}
\begin{definition} Let $\O$ be bounded region in $\R^n$. $\O$ is said to be of {\it A-type}, if there
exists a positive number $A$ such that, for any $x\in\O$ and $0<r<\text{diam}\ \O$,
\[
|\O\cap B(x,r)|\geq Ar^n.
\]
\end{definition}
Now, we state the following Morrey theorem (see, for instance,
\cite{Wulancheng98}):
\begin{theorem}\label{th-morrey}
Assume that $\O$ is of A-type, $u\in W^{1,p}(\O)$, $p>1$ and there
exist constants $K>0,\ 0<\a<1$ such that, for any $B_R$,
\[
\int_{\O\cap B_R}|\g u(x)|^pdx \leq K R^{n-p+\a p}
\]
holds. Then $u\in C^\a(\cl\O)$ and
\[
\text{osc}_{B_R} u\leq CKR^\a,
\]
where $C$ depends on $n, \a, p$ and $A$.
\end{theorem}

\subsection{Uniform Poincar\'e inequality}\

Here, we prove a useful lemma of Poincar\'e type inequality. Assume
that $S$ is a bounded domain in $\R^n$, and there is a constant
$C(n,p,U)$ such that the following classical Poincar\'e inequality
holds:
\[
\bigg(\int_U |u(y)|^p dy\bigg)^{\f 1p} \leq C(n,p,U)\bigg(\int_U
|\g_y u(y)|^p dy\bigg)^{\f 1p},
\]
with $\int_U u(y) dy = 0$, where $C(n,p,U)$ depending only on
$n,p,U$, not on $u$.

Define a class $\s U_K$ by
\[
\ba \s U_K = \big\{\O\big|\ \exists&\text{ an invertible smooth
mapping }
 T: \O\ra U, \text{such\ that }\\
&\|T,\ T^{-1}\|_{C^{2,\alpha}}\leq K<\infty\big\}. \ea
\]
Then, the following useful uniform Poincar\'e type inequality holds:
\begin{proposition}\label{prop-upi}
For any $1\leq p<\infty$, there exists a constant $C(n, p, U, K)$
depending only on  $n, p, U, K$, such that, for any $\O\in \s U_K$,
\[
\int_\O |u(x)|^pdx \leq C(n, p, U, K)\int_\O |\g u(x)|^p dx
\]
or \be\label{eq-88} \|u\|_{L^p(\O)}\leq C(n,p,U, K)\|\g
u\|_{L^p(\O)} \ee holds, provided that $\int_\O u(x) dx = 0$.
\end{proposition}
\proof Set $\a = \f1{|U|}\int_U u\circ T^{-1}(z)dz$ and $J=\f{\p
x}{\p y}$. Since
\[
\int_U u\circ T^{-1}(y) J dy = \int_\O u(x)dx = 0,
\]
one has
\[
\a|\O|=\int_U \left(\a - u\circ T^{-1}(y)\right)J dy.
\]
Then, by the classical Poincar\'{e} inequality for $p=1$ on $U$ , one gets
\[
\ba |\a||\O|&=\bigg|\int_U \left(u\circ T^{-1}(y) - \a\right)J dy\bigg|\\
&\leq \|J\|_{L^{\infty}}\int_U\bigg|u\circ T^{-1}(y) - \a\bigg|dy\\
&\leq \|J\|_{L^{\infty}} C(n,1,U)\int_U \left|\g_y(u\circ T^{-1}(y))\right|dy\\
&\leq \|J\|_{L^{\infty}} C(n,1,U)\left\|\g T^{-1}\right\|_{L^{\infty}}\int_U \left|\g_x u\circ T^{-1}(y)\right|dy\\
&\leq  C(n,1,U)\|\g T^{-1}\|_{L^{\infty}}\|J\|_{L^{\infty}}\bigg(\int_U \left|\g_x u\circ T^{-1}(y)\right|^p dy\bigg)^{\f 1p}|U|^{1-\f 1p}\\
&\leq  C(n,1,U)\|\g T^{-1}\|_{L^{\infty}}\|J\|_{L^{\infty}}\left\|\f
1J\right\|_{L^{\infty}}^{\f 1p}\bigg(\int_\O |\g_x u(x)|^p
dx\bigg)^{\f 1p}|U|^{1-\f 1p}. \ea
\]
Hence
\be\label{eqp1}
\ba
|\a||U|^{\f 1p}&\leq C(n,1,U)\|\g T^{-1}\|_{L^{\infty}}\|J\|_{L^{\infty}}\left\|\f 1J\right\|_{L^{\infty}}^{\f 1p}\|\g_x
u(x)\|_{L^p(\O)}\f{|U|}{|\O|}\\
&\leq C(n,1,U)\|\g T^{-1}\|_{L^{\infty}}\|J\|_{L^{\infty}}\left\|\f
1J\right\|_{L^{\infty}}^{1+\f 1p}\|\g_x u(x)\|_{L^p(\O)}. \ea\ee On
the other hand, by the classical Poincar\'e inequality for $p$ on
$U$, one has
$$
\ba
\bigg(\int_\O |u(x)|^p\f 1J dx\bigg)^{\f 1p} &= \bigg(\int_U |u\circ T^{-1}(y)|^p dy\bigg)^{\f 1p}\\
&\leq \bigg(\int_U |u\circ T^{-1}(y)-\a|^p dy\bigg)^{\f 1p} +|\a||U|^{\f 1p}\\
&\leq C(n,p,U)\bigg(\int_U |\g_y(u\circ T^{-1}(y))|^p dy\bigg)^{\f 1p} +|\a||U|^{\f 1p}\\
&\leq C(n,p,U)\|\g T^{-1}\|_{L^{\infty}}\bigg(\int_\O |\g_x u(x)|^p\f 1J dx\bigg)^{\f 1p} +|\a||U|^{\f 1p}\\
&\leq C(n,p,U)\|\g T^{-1}\|_{L^{\infty}}\left\|\f 1J\right\|_{L^{\infty}}^{\f 1p}\|\g_x u(x)\|_{L^p(\O)} +|\a||U|^{\f 1p}\\
&\leq C(n,p,U)\|\g T^{-1}\|_{L^{\infty}}\|J\|_{L^{\infty}}\left\|\f
1J\right\|_{L^{\infty}}^{1+\f 1p}\|\g_x u(x)\|_{L^p(\O)}
+|\a||U|^{\f 1p},\ea
$$
which, together with (\ref{eqp1}) shows
\[
\ba
\bigg(\int_\O &|u(x)|^p\f 1J dx\bigg)^{\f 1p}
\leq C(n,p,U)\|\g T^{-1}\|_{L^{\infty}}\left\|\f 1J\right\|_{L^{\infty}}^{\f 1p}\|\g_x u(x)\|_{L^p(\O)} +|\a||U|^{\f 1p}\\
&\leq \bigg(C(n,p,U)+C(n,1,U)\bigg)\|\g
T^{-1}\|_{L^{\infty}}\|J\|_{L^{\infty}}\left\|\f
1J\right\|_{L^{\infty}}^{1+\f 1p}\|\g_x u(x)\|_{L^p(\O)}. \ea
\]
Therefore,
\[\ba
\|u(x)\|_{L^p(\O)}&=\bigg(\int_\O |u(x)|^p \f 1J J dx\bigg)^{\f 1p}\\
&\leq \|J\|_{L^{\infty}}^{\f 1p} \bigg(\int_\O |u(x)|^p \f 1J dx\bigg)^{\f 1p}\\
&\leq \bigg(C(n,p,U)+C(n,1,U)\bigg)\|\g
T^{-1}\|_{L^{\infty}}\|J\|_{L^{\infty}}^{1+\f 1p}\left\|\f
1J\right\|_{L^{\infty}}^{1+\f 1p}\|\g_x u(x)\|_{L^p(\O)},\ea\] which
implies the inequality (\ref{eq-88}).
\endproof

\begin{theorem}\label{th-upi}
(Uniform Poincar\'e Inequality) For any $a\in \R$, $1\leq p
<\infty$, one has \be\label{eq-upi} \bigg\|f(x)
-\aint_{\O_{a,a+1}}f(x)dx\bigg\|_{L^p(\O_{a,a+1})} \leq C\|\g
f(x)\|_{L^p(\O_{a,a+1})}. \ee Here
\[
\O_{a,b} = \{x=(x_1,\cdots,x_n)\in\O | a<x_n<b\},
\]
\[
\aint_{\O_{a,a+1}}f(x)dx = \f 1{|\O_{a,a+1}|}
\int_{\O_{a,a+1}}f(x)dx,
\]
$C$ is a positive constant depending only on $n, p, \O$, independent
of $f$, $a$.
\end{theorem}
\proof This can be easily deduced from Proposition \ref{prop-upi}.
\endproof
\subsection{Basic properties of $\O$}

\begin{lemma}\label{L1} Under the assumption (\ref{H1}), for any $0\leq k\leq 2$, if $\fai\in C^{k,\alpha}(\O)$,
$\t\fai=\fai\circ T^{-1}\in C^{k,\alpha}(\C)$
and vise versa, and
\[
\|\fai\|_{C^{k,\alpha}(\cl\O)} \sim \|\t\fai\|_{C^{k,\alpha}(\cl\C)}.
\]
Similar equivalence holds for $H^s$ norms, $s=0,1,2$.
\end{lemma}
\proof The proof follows from simple calculations, and is
omitted.\endproof

According to Lemma \ref{L1}, we may abuse a bit of the notations by
simply denoting $\|\fai\|_{C^{k,\alpha}(\O)}$ and
$\|\t\fai\|_{C^{k,\alpha}(\C)}$ by $\|\fai\|_{C^{k,\alpha}}$ or
$\|\t\fai\|_{C^{k,\alpha}}$, $\|\fai\|_{H^s(\O)}$ and
$\|\t\fai\|_{H^s(\C)}$ by $\|\fai\|_{s}$ or $\|\t\fai\|_s$
respectively.

\begin{lemma}\label{lem-bd}
Assume that $\O$ satisfies (\ref{H1}). Then for any $x_0\in\p\O$,
there exists an invertible $C^{2,\a}$ map $T_{x_0}: U_{x_0}\rightarrow B_{\d_0}: x\mapsto y$
satisfying the following
properties

\begin{subequations}\renewcommand{\theequation}{\theparentequation-\arabic{equation}}\label{eq-bd}
\begin{align}
&\quad T_{x_0}(U_{x_0}\cap\O) = B_{\d_0}^+,\ T_{x_0}(U_{x_0}\cap \p\O) =B_{\d_0}\cap\{y_n=0\},\\
&\smash{%
\raisebox{-0.5\baselineskip}%
{$\displaystyle \left\{
\vphantom{\begin{gathered}x\\x\\x\\x\\x\\x\end{gathered}}%
\right.$}%
} \sigma_{ij}\sigma_{in}(x) = \sigma_{ij}\sigma_{in}(y) =0,\ for\  x\in\p\O\ (\text{i.e. } y_n=0),\quad 1\leq j\leq n-1,\\
&\quad \left\|T_{x_0},T_{x_0}^{-1}\right\|_{C^{2,\a}}\leq K,\\
&\quad |\sigma_{ij}(x)\xi_j|,\ |\sigma_{ij}(y)\xi_j|\sim |\xi|,\quad
\forall x\in U_{x_0}, \forall\ y\in B_{\d_0}^+,\ \forall\ \xi\in
\R^n.
\end{align}
\end{subequations}
where $U_{x_0}$ is a neighbourhood of $x_0$ in $\R^n$, $B_{\d_0}$ is
a ball centered at the origin with radius $\d_0$, $B_{\d_0}^+=
B_{\d_0}\cap\{y_n>0\}$, $\disp \sigma_{ij} = \f{\p y_j}{\p x_i}$,
$\d_0$ and $C$ are positive numbers independent on $x_0\in\p\O$.
\end{lemma}
\proof By assumption (\ref{H1}), $T(x_0) = \t x_0\in\p\C$. Set
\[
V_{x_0} = T^{-1}\left(B_{1/4}(\t x_0)\cap \C\right), \ \
S_{x_0}=\cl{V_{x_0}}\cap\p\O,\ \ \t S_{\t x_0}= \cl{B_{1/4}(\t
x_0)}\cap\p\C,
\]
then
\[
S_{x_0} = T^{-1}(\t S_{\t x_0}).
\]
Suppose that $\t x(y_1,...,y_{n-1}) = \t x(y')$ is the standard surface
parameter of $\t S_{\t x_0}$(and then, of $S_{x_0}$), $\vec N_{x_0}(y')$ is
 the unit inner normal vector on $S_{x_0}$. Let
\[\ld_i(y')=\vec e_i(y')\cdot \vec N_{x_0}(y'),\ \ \ \ 1\leq i\leq n,\]
where $\vec e_i,\ (1\leq i\leq n)$ are the unit coordinate vectors.
Then
\[
\vec N_{x_0}(y') = (\ld_1(y'),...,\ld_n(y')),\ \ \ \ \ld_i(y')\in
C^{1,\a}.
\]

Define $y = T_{x_0}(x)$ by \be\label{eqct} x_i = x_i(y') +y_n^{2-n}
\int_{y_1}^{y_1+y_n}\int_{y_2}^{y_2+y_n}\cdots\int_{y_{n-1}}^{y_{n-1}+y_n}\ld_i(s_1,...,s_{n-1})ds_1ds_2\cdots
ds_{n-1},\ee for $1\leq i\leq n.$ Since $\ca{\ld_i(y')}\leq C\cb{T,\
T^{-1}}$, there exists a $\d_0>0$ independent of $x_0$ such that
$T_{x_0}^{-1}(y)$ is well defined on $B_{\d_0}$ and $\cb{T_{x_0},\
T_{x_0}^{-1}}\leq K$. Then, define $U_{x_0} = T^{-1}(B_{\d_0})$.
Clearly $T_{x_0}$ satisfies (\ref{eq-bd}-1) and (\ref{eq-bd}-3).

Denote the matrix $(\sigma_{ij}(x))$ by $A(x)$. For any $\xi\in \R^n$,
\[
|A(x)\xi|\leq |A(x)||\xi|\leq C|\xi|,\quad |A^{-1}(x)\xi|\leq |A^{-1}(x)||\xi|\leq C|\xi|.
\]
Then (\ref{eq-bd}-4) follows immediately.

To prove (\ref{eq-bd}-2), we differentiate (\ref{eqct}) with respect
to $y_j$ and note that $y = (y',0)$ for $x\in\p\O$,
\[
\ba \f{\p x_i}{\p y_j}(y',0) &= \f{\p x_i(y')}{\p y_j} +
\ld_i(y')\d_{jn}. \ea
\]
Since
\[
A^{-1}(y',0)= (\sigma_{ij})^{-1}(y',0) = \bigg(\f{\p y_{i}}{\p x_j}(y',0)\bigg)^{-1} = \bigg(\f{\p x_i}{\p y_j}(y',0)\bigg),
\]
hence
\be\label{eq1}
\bigg(\sigma_{ji}\cdot \f{\p x_i}{\p y_n}\bigg)(y',0) = \sigma_{ji}(y',0)\cdot\ld_i(y',0) = 0,\quad 1\leq j\leq n-1.
\ee

On the other hand, since
\[
\sigma_{ni}(y',0)\cdot\f{\p x_i}{\p y_j}(y',0)= 0,\quad 1\leq j\leq n-1,
\]
$(\sigma_{n1}(y',0),...,\sigma_{nn}(y',0))$ is the normal direction
of $S_{x_0}$ at $x\in S_{x_0}$, that is, $(\sigma_{n1}(y',0),...,\sigma_{nn}(y',0))$ is parallel
to the inner normal $\vec N_{x_0}(x) = (\ld_1(y',0),...,\ld_n(y',0))$. Comparing with (\ref{eq1}) yields
\[
\sigma_{ji}(y',0)\sigma_{ni}(y',0) = 0,\quad 1\leq j\leq n-1.
\]
\endproof

\begin{remark} Hypothesis (\ref{H1}) is stronger than the $C^{2,\a}$-regularity hypothesis on $\O$.
If one only assumes that $\O\in C^{2,\a}$, then $C$ and $\d_0$ in lemma
\ref{lem-bd}, in general, may depend on $x_0\in\p\O$.
\end{remark}

\begin{lemma}\label{lem-in}
There exists a $\d_1>0$ such that
\[
\ba &\d_0\sim \d_1, &\O = \bigg(\bigcup_{x_0\in\p\O}
T_{x_0}^{-1}\left(B_{\f{\d_0}{2}}^+\right)\bigg)\bigcup\bigg(\bigcup_{B_{2\d_1}\subset
\O}B_{\d_1}\bigg). \ea
\]
where $T_{x_0}$ and $\d_0$ are the same as in Lemma \ref{lem-bd}.
\end{lemma}
\proof Since $\cb{T,\ T^{-1}}\leq C$, $|x_1-x_2|\sim |T(x_1) - T(x_2)|$, there exists a constant $\d\sim\d_0$ such that
\[
B_{\d}(x_0)\cap\O \subset
T_{x_0}^{-1}\left(B_{\f{\d_0}{2}}^+\right),\quad \forall\
x_0\in\p\O.
\]
Taking $\d_1 = \f\d 2$ yields the Lemma.
\endproof

\section{The existence of subsonic flow for small incoming mass flux}
There are two major obstacles to solve the Problem \ref{op}. First, the
ellipticity of the equation (\ref{nozzle}) is not guaranteed
beforehand, since there is no a priori $L^\infty$ bound for
$\g\fai$, the gradient of the solution to the Problem \ref{op}. Second, the
nozzle region is unbounded, and can not be transformed
to a bounded domain by Kelvin-like transformations. In order to
overcome these difficulties, we first truncate the coefficients of
the equation in (\ref{nozzle}) to ensure the strong ellipticity, and
then, truncate the domain $\O$ to a series of bounded domains
$\O_L$, with additional boundary conditions. Therefore, to solve the
Problem \ref{op} becomes to study a series of approximate strong elliptic
problems in bounded domains and their uniform estimates, which
ensure to pass the limit of the approximate solutions to the
Problem \ref{op}.

\subsection{A subsonic truncation and approximate solutions}
\subsubsection {A subsonic truncation.}
By normalizing the equation if necessary \cite{BERS},
\cite{XX1}, one can assume that the critical sound speed of the flow
is one. Thus, the density-speed relation (\ref{strongber}), $\rho =
\rho(q^2)$, is positive, sufficiently smooth and nonincreasing in  $q=|\g\fai|\in [0, 1]$. However, the potential equation is not uniformly elliptic as $q$ approaches to 1. To guaranteed the
uniformly ellipticity, we truncate the coefficients as follows.

Define two functions $\H(s^2)$ and $F(q^2)$ as follows
\begin{equation}\label{eqH}
\H(s^2) = \left\{
\begin{array}{ll}
 \rho(s^2),&\quad\text{ if } s^2 < 1-2\t\delta_0,\\
 \text{monotone and smooth}, &\ \ \ \text{ if }1-2\tilde\delta_0\leq s^2\leq 1-\tilde\delta_0,\\
 \rho(1-\t\delta_0), &\quad\text{ if } s^2 > 1-\t\delta_0,
\end{array}
\right.
\end{equation}
and
\begin{equation}\label{eqF}
 F(q^2) = \f12\int_0^{q^2} \H(s^2) ds^2,
\end{equation}
where $\t\d_0>0$. Moreover,  $\H(s^2)$ is a smooth non-increasing
functions and $F(q^2)$ is a smooth increasing function.
Set
\[
 a_{ij}(\g\fai) = \H(|\g\fai|^2)\d_{ij} + 2\H'(|\g\fai|^2)\p_i\fai\p_j\fai.
\]
It is easy to check the following facts, \be\label{fact1}
 F\left(q^2\right) \sim q^2,\ \ \  \f 1{C(\t\d_0)} <\H(s^2),\ \ \ \H(s^2) +2\H'(s^2) s^2 < C(\t\d_0),\ee
and there exist two positive constants $\ld$ and $\Ld$, such that
\be\label{fact2}
  \ld|\xi|^2<a_{ij}(\g\fai)\xi_i\xi_j<\Ld|\xi|^2,
\ee where $C(\t\d_0)$, $\ld$ and $\Ld$ depend only on the subsonic
truncation parameter $\t\d_0$. Note that a solution of the potential
equation derived from the new density-speed relation $\H(q^2)$ is
also a solution of the actual potential equation provided that
$|\g\fai|^2 \leq 1-2\t\d_0$. Therefore, in the end of this section,
we will show that the solution of the truncated problem satisfies
$|\g\fai|^2 \leq 1-2\t\d_0$, as long as the incoming mass flux $m_0$
is suitable small. Consequently, the subsonic truncation can be
removed.

\subsubsection{Domain truncation.}
Our strategy to deal with the unbounded domain here is to construct
a series of truncated problems to approximate the Problem \ref{op}
with subsonic truncation.

Let $L>0$ be sufficiently large. Define
\[
 \Omega_L=\big\{x\in \Omega\ \big|\ |x_n| < L\big\},\ \ \ S^\pm_L =\Omega\cap\{x_n=\pm L\},\ \ \ S_L = \sla\cup\slb.
\]
Consider the following truncated problem with $m_0>0$.
\begin{problem}\label{ptnozzle} Find a $\fai$ such that,
\begin{equation}\label{eqTD}
\left\{\begin{array}{ll}
\text{div}(\H(|\g\fai|^2)\g\fai) = 0,\ \ \ \ & x\in\ol,\\
\f{\p \fai}{\p \vec n} = 0,\ \ \ \ &\p\Omega\cap\p\Omega_L,\\
\H(|\g\fai|^2)\frac{\p\fai}{\p x_n} =\f{m_0}{|S_L^+|},& \text{on}\ \ S_L^+\\
\fai = 0,&\text{on} \ \ S_L^-.
\end{array}\right.
\end{equation}
\end{problem}
\noindent The additional boundary condition on $S^+_L$ implies the
mass flux of the flow remains $m_0$.

 Clearly, the truncated problem
\ref{ptnozzle} is a strong quasilinear elliptic problem in a bounded
domain. From now on, instead of the original Problem \ref{op}, we consider a
series of the truncated Problem \ref{ptnozzle}  for any fixed
sufficiently large $L$. With some uniform estimates of the
approximate solutions, we can conclude that the solution of the
truncated problem \ref{ptnozzle} converges to the original Problem \ref{op}.

\subsection{Truncated variational problem} \ \\
\par
In this subsection, we solve the truncated problem \ref{ptnozzle} by
a variational method.
Define
\[
 H_L = \left\{\fai\in H^1(\ol) :\ \fai\big|_{S_L^-} = 0\right\}.
\]
Then, $H_L$ is a Hilbert space under $H^1$-norm. The additional
boundary condition on $S^-_{L}$ is understood in the sense of
traces. Define a functional $J(\psi)$ on $H_L$ as
\[
 J(\psi) = \int_{\ol}F(|\g\psi|^2) dx - \f{m_0}{|S_L^+|}\int_{S_L^+}\psi
 dx',
\]
where $F(q^2)$ is defined by (\ref{eqF}) and $x'=(x_1, x_2, \ldots, x_{n-1})$. The existence of solution to
problem \ref{ptnozzle} is equivalent to the following variational
problem:

\begin{problem}\label{P1}
Find a minimizer $\fai\in H_L$ such that
\[
 J(\fai) = \min\limits_{\psi\in H_L} J(\psi).
\]
\end{problem}

\begin{theorem}\label{th-ae} %%average estimate
 Problem \ref{P1} has a nonnegative minimizer $\fai\in H_L$. Moreover,
 \begin{equation}\label{L}
\frac{1}{|\ol|}\int_{\ol}|\g\fai|^2dx \leq Cm_0^2,
\end{equation}
where the constant $C$ does not depend on $L$.
\end{theorem}

\proof Step 1. $J(\psi)$ is coercive on $H_L$. In fact, by Lemma
\ref{L1}, for any $\psi\in H_L$,
\be\label{n11}
 \begin{aligned}
 \bigg|\int_{S_L^+}\psi dx'\bigg| &\leq C\bigg|\int_{B(0,1)} \t\psi
 dy'\bigg|\leq C\bigg|\int_{B(0,1)} \int_{-L}^L\t\p_{n}\t\psi dy_ndy'\bigg|\\
 &\leq C\int_{\C_L}|\t\g\t\psi|dy\leq C\int_{\ol}|\g\psi|dx\\\
 &\leq C|\ol|^{\f 12}\lb{\g\psi}
 \end{aligned}
\ee
Therefore, applying (\ref{n11}) and Cauchy inequality yields
\[
 \begin{aligned}
 J(\psi) &= \int_{\ol}F(|\g\psi|^2)dx - \f{m_0}{|S_L^+|}\int_{S_L^+}\psi dx'\\
 &\geq
 \ld \int_{\ol}|\g\psi|^2dx - C(m_0, |S^+_L|,|\ol|) \lb{\g\psi}\\
 &\geq \f \ld{2}\lb{\g\psi}^2 - \f 1\ld C(m_0, |S^+_L|,|\ol|),
 \end{aligned}
\]
which implies $J(\psi)$ is coercive.

Step 2. The existence of the minimizer $\fai\in H_L$. Since
$J(\psi)$ is coercive in $H_L$, there is a minimizer sequence
$\{\fai_n\}\subset H_L$ such that
\[
J(\fai_n)\ra \alpha = \inf\limits_{\psi\in H_L}J(\psi)>-\infty.
\]
Then,
\[\ba
 \lb{\g\fai_n}^2 &\leq \f 2 \ld J(\fai_n) + \f 2{\ld^2}C(m_0,
 |S^+_L|,|\ol|)\\
 &\leq \f 2\ld J(0) + \f 2{\ld^2}C(m_0, |S^+_L|,|\ol|)\\
 &= \f 1{\ld^2}C(m_0, |S^+_L|,|\ol|).
\ea\]
Therefore, there exists a subsequence, denoted by $\{\fai_n\}$ converges weakly to some $\fai\in H_L$ and
\[
\lb{\g\fai}^2 \leq \f 1{\ld^2}C(m_0, |S^+_L|,|\ol|).
\]
By Fatou's Lemma, it is easy to check that
\be\label{n12}
 \int_{\ol} F(|\g\fai|^2)dx \leq \liminf_{n\ra\infty}\int_\ol
 F(|\g\fai_n|^2)dx.
\ee
On the other hand,
\[
 \begin{aligned}
\int_\slb(\fai_n-\fai)^2 dx' &\leq C\int_{B(0,1)}(\t\fai_n -
\t\fai)^2 dy'\\
&\leq C\bigg|\int_{B(0,1)}\int_{-L}^{L}(\t\fai_n-\t\fai)\t\p_{n}(\t\fai_n -\t\fai)dy_ndy'\bigg|\\
&\leq C\int_{\C}|\t\fai_n-\t\fai||\t\g\t\fai_n-\t\g\t\fai|dy\\
&\leq C\int_\ol|\fai_n-\fai||\g\fai_n-\g\fai|dx\\
&\leq C\bigg(\int_\ol|\fai_n-\fai|^2dx\bigg)^{\f
12}\bigg(\int_\ol|\g\fai_n-\g\fai|^2dx\bigg)^{\f 12} \ra 0,
\end{aligned}
\]
as $n\ra \infty$. Then, \be\label{n13}
 \int_\slb|\fai_n-\fai|dx' \ra 0, \ \ \ \ \text{ as } n\ra \infty.
\ee Therefore, it follows from (\ref{n12}) and (\ref{n13}) that
\[
 J(\fai)\leq\liminf_{n\ra\infty} J(\fai_n) = \alpha.
\]
i.e.
\[
 J(\fai) = \min\limits_{\psi\in H_L} J(\psi) = \alpha.
\]
Step 3. $\fai^+ = \max\{\fai, 0\}$ is a nonnegative minimizer in
$H_L$. Indeed, since $\fai\in H_L$, $\fai^+ \in H_L$, and
\[
 \begin{aligned}
& |\g\fai^+|^2 \leq |\g\fai|^2,\quad \ \ F(|\g\fai^+|^2)\leq F(|\g\fai|^2),\\
&\f{m_0}{|\slb|}\int_\slb \fai ^+ dx' \geq \f{m_0}{|\slb|}\int_\slb
\fai  dx'.
 \end{aligned}
\]
Hence,
\[
 J(\fai^+)\leq J(\fai).
\]
Since $\fai$ is a minimizer, $J(\fai^+) = J(\fai)$, which implies
that $\fai^+ \geq 0$ is also a minimizer.

\noindent
Step 4. By direct computations,
\[
 \begin{aligned}
 \int_\ol F(|\g\fai|^2)dx &= J(\fai) + \f{m_0}{|\slb|}\int_\slb\fai dx'\leq J(0) + \f{m_0}{|\slb|}\int_\slb\fai dx'\\
 &\leq C\f{m_0}{|\slb|}|\ol|^{\f 12}\lb{\g\fai}.
 \end{aligned}
\]
It follows from (\ref{fact1}) and (\ref{fact2}) that
\[
 \lb{\g\fai}^2\leq \f 1\ld \int_\ol F(|\g\fai|^2)dx \leq C\f 1\ld\f{m_0}{|\slb|}|\ol|^{\f
 12}\lb{\g\fai}.
\]
That is
\[
 \lb{\g\fai}^2 \leq C\f {m_0^2}{\ld^2 |\slb|^2}|\ol|,
\]
i.e.
\[
 \f 1{|\ol|}\int_\ol |\g\fai|^2dx \leq C\f{m_0^2}{\ld^2|\slb|^2} \leq
 C\f{m_0^2}{\ld^2S_{min}^2},
\]
where $S_{min}$ denotes the minimal of $|\slb|$.
\endproof
\begin{remark}
The estimate (\ref{L}) is the key estimate for the existence of the
classical solution to Problem \ref{ptnozzle}. Indeed, the potential
$\fai$ is essentially unbounded, one can not expect to get uniform
bounds on $\|\g\fai\|_{L^\infty}$ through $\|\fai\|_{L^\infty}$ as in
the standard elliptic theory.
\end{remark}

\begin{proposition}\label{th-ws} %%weak solution
$\fai\in H_L$ is a weak solution to the equations in (\ref{eqTD}) in the following sense:
\begin{equation}\label{eq-W}
\int_\ol \H(|\g\fai|^2)\g\fai\cdot\g\psi dx -
\f{m_0}{|\slb|}\int_\slb \psi dx' = 0,\quad\quad \forall\ \psi\in
H_L
\end{equation}
\end{proposition}
\proof This is a standard variation problem. In fact, for any
$t\in \R,\ t>0$ and any $\psi\in H_L$, $\fai+t\psi\in H_L$. Then,
\be\label{n14}
 \begin{aligned} %%
 0 &\leq J(\fai + t\psi) - J(\fai)
  = \int_\ol F(|\g\fai + t\g\psi|^2)-F(|\g\fai|^2)dx - \f{m_0t}{|\slb|}\int_\slb\psi dx'.
 \end{aligned}
\ee Mean value theorem yields that
\be\label{n15}
\begin{aligned}
& \ \ \ \ \int_\ol F(|\g\fai + t\g\psi|^2)-F(|\g\fai|^2)dx\\
&= \int_\ol \int_0^1 F'(\theta|\g\fai + t\g\psi|^2 + |\g\fai|^2(1-\theta))d\theta (|\g\fai+t\g\psi|^2-|\g\fai|^2)dx\\
&= \int_\ol \int_0^1 F'(|\g\fai|^2 +\theta(t^2|\g\psi|^2+2t\g\psi\cdot\g\fai))d\theta (t^2|\g\psi|^2+2t\g\psi\cdot\g\fai)dx.
\end{aligned}
\ee Since $|F'(\cdot)|\leq C$, $\g\fai,\ \g\psi\in L^2(\ol)$,
substituting (\ref{n15}) into (\ref{n14}) shows that
\[
\begin{aligned}
0&\leq \liminf_{t\ra 0^+}\f 1t(J(\fai + t\psi) - J(\fai))\\
&= \liminf_{t\ra 0^+}\int_\ol \int_0^1 F'(|\g\fai|^2 +\theta(t^2|\g\psi|^2+2t\g\psi\cdot\g\fai))d\theta (2\g\psi\cdot\g\fai)dx\\
&~~~~~ -\int_\slb\f{m_0}{|\slb|}\psi dx'\\
&=\int_\ol \int_0^1 F'(|\g\fai|^2)d\theta (2\g\psi\cdot\g\fai))dx -\f{m_0}{|\slb|}\int_\slb\psi dx'\\
 &~~~~~~~~\text{(by Lebesgue's theorem)}\\
&=\int_\ol \H(|\g\fai|^2)\g\fai\cdot\g\psi dx
-\f{m_0}{|\slb|}\int_\slb\psi dx'.
\end{aligned}
\]
Therefore, for any $\psi\in H_L$,
\[
 \int_\ol \H(|\g\fai|^2)\g\fai\cdot\g\psi dx -\f{m_0}{|\slb|}\int_\slb\psi dx' = 0.
\]
\endproof

\subsection{$H^2$ regularity of the weak solution}

We are now ready to improve the regularity of the minimizer $\fai$. Indeed, one has
\begin{proposition}\label{th-h2} %%h2-regularity
$\fai\in H^2\left(\O_{L/2}\right)$. Moreover, \be\label{eq-wb}
\t\p_n\t\fai(y)\bigg|_{y_n = 0} = \f{\p\t\fai(y)}{\p y_n}\bigg|_{y_n
= 0} = 0. \ee
\end{proposition}
To prove this, one needs the following estimates in Lemma
\ref{lem-H2i}, \ref{lem-H2b}.
\begin{lemma}\label{lem-H2i} %% inner H2 regularity
(Interior Estimate) For any $B_{2R}(x_0)\subset \ol,\ R\leq \d_1$,
\[
\g^2\fai\in L^2(B_R).
\]
Here $\d_1$ is the same number in the Lemma \ref{lem-in}.
\end{lemma}
\proof For any $B_{2R}(x_0)\subset \ol$, $v \in H_0^1\left(B_{\f32
R}\right)$, $h<\f12R$, one has \be\label{1eq62} \ba
0&=\int_{B_{2R}}\H(|\g\fai|^2)\g\fai\cdot\g(\d_{-h}v)dx =
-\int_{B_{2R}}\d_{h}(\H(|\g\fai|^2)\g\fai)\cdot\g v dx, \ea \ee
where $\d_h v(x)\eqdef \f 1h(v(x+h\vec {e}_k)-v(x))$ is the $k$-th
difference quotient, $k=1,2,\cdots,n$.

Set
$$
\t q = t\g\fai^h + (1-t)\g\fai,\quad\fai^h(x) = \fai(x+h\vec e_k),
$$
\[
a_{ij}(\t q, t) = \H(\t q^2)\d_{ij} +2\H'(\t q^2)\t q_i\t q_j,\quad
a_{ij}=a_{ij}(\t q) = \int_0^1 a_{ij}(\t q, t)dt.
\]
Then, direct calculations give \be\label{n17}
\d_h(\H(|\g\fai|^2)\g\fai) =\int_0^1a_{ij}(\t q, t)dt\p_j(\d_h\fai)
= a_{ij}\p_j(\d_h\fai). \ee Therefore, substituting (\ref{n17}) into
(\ref{1eq62}), one has \be\label{eq-93} \int_{B_{2R}}
a_{ij}\p_j(\d_h\fai)\p_jv dx = 0. \ee

Take $v = \eta^2\d_h \fai$ in (\ref{eq-93}), where $\eta\in
C_0^\infty\left(B_{\f32R}\right)$, $\eta \equiv 1$ in $B_R$,
$|D\eta|\leq \f C R$. Then,
\be\label{n18} \ba
0&=\int_{B_{2R}}a_{ij}\p_j(\d_h\fai)\p_j(\eta^2\d_h\fai)dx\\
 &=\int_{B_{2R}}\eta^2a_{ij}\p_j(\d_h\fai)\p_j(\d_h\fai)dx +
 2\int_{B_{2R}}a_{ij}\p_j(\d_h\fai)\eta\p_j\eta\d_h\fai dx.
\ea \ee It follows from H\"{o}lder inequality and (\ref{n18}) that
\[
\ba & \ \ \ \int_{B_{2R}}\eta^2a_{ij}\p_j(\d_h\fai)\p_j(\d_h\fai)dx
=-2\int_{B_{2R}}a_{ij}\p_j(\d_h\fai)\eta\p_j\eta\d_h\fai dx\\
&\leq
2\bigg(\int_{B_{2R}}\eta^2a_{ij}\p_i(\d_h\fai)\p_j(\d_h\fai)dx\bigg)^{\f12}
\bigg(\int_{B_{2R}}a_{ij}\p_i\eta\p_j\eta(\d_h\fai)^2dx\bigg)^{\f12},
\ea
\]
namely,
\[
\int_{B_{2R}}\eta^2a_{ij}\p_i(\d_h\fai)\p_j(\d_h\fai)dx \leq
4\int_{B_{2R}}a_{ij}\p_i\eta\p_j\eta(\d_h\fai)^2dx.
\]
Consequently, by the strong ellipticity of $a_{ij}$, one gets
\[
\ba
\ld\int_{B_{R}}|\g(\d_h\fai)|^2dx &\leq \int_{B_{2R}}a_{ij}\p_i(\d_h\fai)\p_j(\d_h\fai)dx\\
&\leq 4\int_{B_{2R}}a_{ij}\p_i\eta\p_j\eta(\d_h\fai)^2dx\\
&\leq C\f{\Ld}{R^2}\int_{B_{\frac{3}2R}}(\d_h\fai)^2dx\\
&\leq C\f{\Ld}{R^2} \int_{B_{\f{3}2R}}|\g\fai|^2dx, \ea
\]
and then \be \label{n19} \int_{B_{R}}|\g(\d_h\fai)|^2dx \leq
C\f{\Ld}{\ld R^2}\int_{B_{2R}}|\g\fai|^2dx,\quad \forall\ h<R. \ee
Therefore, according to (\ref{n19}) and $H^1$ regularity of
minimizer $\fai$, we can conclude that $\g^2\fai\in L^2(B_R)$.
\endproof

Next, we derive the boundary estimate of the minimizer $\fai$.
\begin{lemma}\label{lem-H2b} %%H2 boundary regularity
(Boundary Estimate) For any $x_0\in\p\O_{L/2}$, \be\label{1eq-be}
\g^2\fai\in L^2\left(B_{\f{\d_0}2}(x_0)\cap\ol\right). \ee
\end{lemma}
\proof \noindent Set $U_{x_0,\delta_0}= B_{\delta_0}(x_0)\cap \O_L$,
and  $$T_{x_0}: U_{x_0,\delta_0}\ra B_{\d_0}^+: x\mapsto y,~~~
y=T_{x_0}(x),~~  \sigma_{ij}(y) = \f{\p y_j}{\p x_i}(y),~~ J(y) =
\f{\p x}{\p y}.$$ For simplification, we write $U_{x_0,\delta_0}$
and $B_{\d_0}^+$ as $U$ and $B^+$ respectively in the remaining of
the proof.

Then for any $\psi\in H_0^1(U)$, $\t\psi=\psi\circ T^{-1}_{x_0}$,
\[
0=\int_U \H(|\g\fai|^2)\g\fai\cdot\g\psi dx =
\int_{B^+}\H(|\sigma_{\a\b}\t\p_\b\t\fai|^2)\sigma_{ij}\t\p_j\t\fai\sigma_{il}\t\p_l\t\psi
J dy,
\]
where $\t\fai=\fai\circ T^{-1}_{x_0}$. Taking $\t\psi$ as the $k$-th
difference quotient $$\t\d_{-h}\t\psi\eqdef \f 1h\left(\t\psi(y) -
\t\psi(y-h\vec e_k)\right)~~~~~ \text{for}~~~k=1,2,\cdots,n-1,$$ we
may get from the property and the "integrate by parts" formula for
difference quotient that for suitable small $h>0$
\[
\ba
0 &=\int_{B^+}\t\d_h\left(\H(|\sigma_{\a\b}\t\p_\b\t\fai|^2)\sigma_{ij}\sigma_{il}J\t\p_j\t\fai\right)\t\p_l\t\psi dy\\
&=\int_{B^+}\t\d_h\left(\H(|\sigma_{\a\b}\t\p_\b\t\fai|^2)\sigma_{ij}\t\p_j\t\fai\right)\sigma_{il}\t\p_l\t\psi
J dy\\
&~~~+\int_{B^+}\left(\H(|\sigma_{\a\b}\t\p_\b\t\fai|^2)\sigma_{ij}\t\p_j\t\fai\right)^h\t\d_h\left(\sigma_{il}J\right)\t\p_l\t\psi dy\\
&= I + II. \ea
\]
Set
$$\t A_{ij} = \int_0^1 \t A_{ij}(t)dt,\quad \t A_{ij}(t) = \H(|q(t)|^2)\d_{ij} + 2\H'(|q(t)|^2)q_i(t)q_j(t),$$
$$q(t) = tq^h(y) +(1-t)q(y)=(q_1(t),q_2(t),\cdots,q_n(t)),$$$$\quad q^h(y) = \sigma_{\a\b}(y+h\vec
e_k)\t\p_\b\t\fai(y+h\vec e_k),~~~~
q(y) =
\sigma_{\a\b}(y)\t\p_\b\t\fai(y).$$
Now, the term $I$ can be rewritten as
\[
\ba I
&=\int_{B^+}\t\d_h\left(\H(|\sigma_{\a\b}\t\p_\b\t\fai|^2)\sigma_{ij}\t\p_j\t\fai\right)\sigma_{il}\t\p_l\t\psi
Jdy\\
&=\int_{B^+}\t A_{ij}\f 1h\left(q^h_j-q_j\right)\sigma_{il} \t\p_l\t\psi J dy\\
&=\int_{B^+}\t A_{ij}\sigma_{js}\t\p_s(\t\d_h\t\fai)\sigma_{il}
\t\p_l\t\psi Jdy
+\int_{B^+}\t A_{ij}\left(\t\p_s\t\fai\right)^h\t\d_h (\sigma_{js})\sigma_{il} \t\p_l\t\psi Jdy\\
&=I_1 + I_2. \ea
\]
Set
$$
\t\psi =  \t\eta^2 \t u_h,\quad \t u_h = \t\d_h\t\fai,$$
\[\t\eta\in C^\infty_0(B^+),\quad \t\eta\equiv 1\ \ \text{in}\ \ \t
B^+=T_{x_0}\left(U_{x_0,\frac{\d_0}2}\cap \O_L\right),\quad
|\t\g\t\eta| \leq 2
\]
and
\[
\psi = \t\psi \circ T_{x_0}= \eta^2 u_h,\quad u_h = \t u_h\circ
T_{x_0},
\]
\[
\eta=\t\eta\circ T_{x_0}\in C^\infty_0(U\cap\O_L), \quad \eta\equiv
1\ \ \text{in}\ \ U_{x_0,\frac{\d_0}{2}}\cap \O_L,\quad |\g\eta|\leq
2.
\]
Then
\[
\ba I_1 &=\int_{B^+}\t A_{ij}\sigma_{js}\t\p_s \t u_h \sigma_{il}\t\p_l\t\psi Jdy\\
&=\int_U \t A_{ij} \p_i\psi \p_j u_h dx\\
&=\int_U \t A_{ij}\eta^2\p_i u_h \p_j u_h  dx + 2\int_U \t A_{ij}\eta u_h\p_j u_h  \p_i\eta  dx\\
&= I_{11} + I_{12}, \ea
\]
and
\[
\ba I_2 &=\int_{B^+}\t A_{ij}(\t\p_s\t\fai)^h\t\d_h
(\sigma_{js})\sigma_{il} \t\p_l \t u_h\t\eta^2 J dy
+2\int_{B^+}\t A_{ij}(\t\p_s\t\fai)^h\t\d_h (\sigma_{js})\sigma_{il} \t u_h\t\eta\t\p_l\t\eta J dy\\
&= I_{21} + I_{22}. \ea
\]
Due to the strong ellipticity, \be\label{n20} I_{11}=\int_{U}\t
A_{ij}\eta^2\p_iu_h\p_ju_hdx \geq \lambda \|\eta\g u_h\|_{L^2(U)}^2.
\ee To estimate the term $I_{11}$, we will deal with $I_{12}$,
$I_{21}$, $I_{22}$ and $II$ first.

By H\"older inequality and the strong ellipticity of $\t A_{ij}$, we
have
\[
\ba |I_{12}| &=\bigg|2\int_U \t A_{ij}\eta u_h\p_j u_h  \p_i\eta
dx\bigg|\\
&\leq 2\bigg(\int_U \t A_{ij} \eta^2\p_i u_h \p_j u_hdx\bigg)^{\f12}
\bigg(\int_U \t A_{ij}u_h^2\p_j\eta\p_i\eta dx\bigg)^{\f 12}\\
&\leq C I_{11}^{\f 12} \bigg(\Ld \|\g\eta\|_{L^\infty}^2\cdot
\|u_h\|^2_{L^2(U)}\bigg)^{\f 12}\\
&\leq C I_{11}^{\f 12} \bigg(\Ld
\|\g\fai\|^2_{L^2(U)}\bigg)^{\f 12}\\
&\leq \f14 I_{11} + C \Ld  \|\g\fai\|^2_{L^2(U)}, \ea
\]
and
\[
\ba
|I_{21}|&=\bigg|\int_{B^+}\t A_{ij}(\t\p_s\t\fai)^h\t\d_h (\sigma_{js})\sigma_{il} \t\p_l \t u_h\t\eta^2 J dy\bigg|\\
&\leq\bigg(\int_{B^+}\t A_{ij}\sigma_{il}\t\p_l\t u_h
\sigma_{js}\t\p_s\t u_h\t\eta^2 Jdy\bigg)^{\f 12}
\bigg(\int_{B^+}\t A_{ij}\t\eta^2 \t\d_h(\sigma_{js})(\t\p_s\t\fai)^h\t\d_h (\sigma_{il})(\t\p_l\t\fai)^h J dy\bigg)^{\f 12}\\
&\leq C I_{11}^{\f 12} \left(\Ld |K|^2
\|\g\fai\|^2_{L^2(U)}\right)^{\f12}\\
&\leq \f14 I_{11} + C \left(\Ld |K|^2 \|\g\fai\|^2_{L^2(U)}\right),
\ea
\]
where $K$ is $C^{2,\a}$ norm of the boundary (See assumption
(\ref{H1})).
\[
I_{22} =2\int_{B^+}\t A_{ij}(\t\p_s\t\fai)^h\t\d_h
(\sigma_{js})\sigma_{il} \t u_h\t\eta\t\p_l\t\eta J dy \leq C \Ld
K\|\g\fai\|^2_{L^2(U)}.
\]
Next, we estimate $II$.
\[
\ba
II &= \int_{B^+}\left(\H(|\sigma_{\a\b}\t\p_\b\t\fai|^2)\sigma_{ij}\t\p_j\t\fai\right)^h\t\d_h(\sigma_{il}J)\t\p_l(\eta^2 \t u_h) dy\\
&=\int_{B^+}\t\eta^2\left(\H(|\sigma_{\a\b}\t\p_\b\t\fai|^2)\sigma_{ij}\t\p_j\t\fai\right)^h\t\d_h(\sigma_{il}J)\t\p_l \t u_h dy\\
&\ \ \ +2\int_{B^+}\t\eta\left(\H(|\sigma_{\a\b}\t\p_\b\t\fai|^2)\sigma_{ij}\t\p_j\t\fai\right)^h\t\d_h(\sigma_{il}J)\t u_h\t\p_l\t\eta dy\\
&=II_1+II_2. \ea
\]
Then direct computations yield that
\[
\ba
II_1 &=\int_{B^+}\t\eta^2\left(\H(|\sigma_{\a\b}\t\p_\b\t\fai|^2)\sigma_{ij}\t\p_j\t\fai\right)^h\t\d_h(\sigma_{il}J)\t\p_l \t u_h dy\\
&\leq C\Ld K\int_{B^+}\t\eta^2|\t\g\t\fai||\t\g \t u_h|dy\\
&\leq C\Ld K\int_{U}\eta^2|\g\fai||\g u_h|dx\\
&\leq \f{\ld}4\|\eta\g u_h\|^2_{L^2(U)} + \f{C}{\ld}\Ld^2
K^2\|\g\fai\|^2_{L^2(U)}, \ea
\]
and
\[
\ba II_2&\leq 2\Ld K\int_{B^+}|\t\g\t\fai||\t\g\t\eta||\t u_h|dy
\leq C\Ld K \|\g\fai\|^2_{L^2(U)}. \ea
\]
Therefore, noticing that
$$
0=I_1+I_2+II_1+II_2=I_{11}+I_{12}+I_{21}+I_{22}+II_1+II_2,
$$
and applying the estimates as above, we get
\[\ba
I_{11}&\leq C\Ld (K^2+1)\|\g\fai\|^2_{L^2(U)}
+\f \ld 2 \|\eta\g u_h\|^2_{L^2(U)} + \f{C}{\ld}\Ld^2 K^2\|\g\fai\|^2_{L^2(U)}\\
&\leq \f \ld 2 \|\eta\g u_h\|^2_{L^2(U)} + C
(K^2+1)
\|\g\fai\|^2_{L^2(U)}. \ea\] Then, combining with (\ref{n20}), we
obtain the gradient estimates for the $k$th difference quotient $\t
u_h$ ($k=1,2,\cdots,n-1$),
$$
\|\t\g\t u_h\|^2_{L^2(\t B^+)} \leq C (K^2+1)
 \|\g\fai\|^2_{L^2(U)}.
$$
Furthermore, the following derivatives estimates hold,
\be\label{n21} \sum_{k=1}^{n-1}\|\t\g(\t D_k\t\fai)\|^2_{L^2\left(\t
B^+\right)}\leq C (K^2+1)\|\g\fai\|^2_{L^2(U)}. \ee For the $\t
D^2_{nn} \t\fai$, by the potential equation, and the estimates for
$\t D^2_{kj},\ 1\leq k\leq n-1,\ 1\leq j\leq n$, \be\label{n22} \|\t
D^2_{nn} \t\fai\|^2_{L^2\left(\t B^+\right)}\leq C (K^2+1)
\|\g\fai\|^2_{L^2(U)}. \ee Combining the estimates (\ref{n21}),
(\ref{n22}) and the $H^1$-estimate (\ref{L}) yields (\ref{1eq-be}).
\endproof

\noindent{\sl Proof of Proposition \ref{th-h2}}: It follows from
Lemma \ref{lem-H2i}, Lemma \ref{lem-H2b} and a finite cover
argument. {\hfill$\Box$}

\subsection{Local average estimate} Set
\[
\O_{x_0,r} = \O\cap\{x=(x',x_n):\ |x_n-x_{0,n}|<r\}, \text{ where } x_0=(x_0', x_{0,n})\in \O.
\]
\begin{proposition}\label{th-lae} %%local average estimate
{(Local average estimate).} For any $x_0\in\O$ with
$|x_{0,n}|<\f 12 L$, one has
\begin{equation}\label{eq-LAE} %%Local Average Estimate
 \f 1{|\O_{x_0,1}|}\int_{\O_{x_0,1}}|\g\fai|^2 dx \leq Cm_0^2,
\end{equation}
where $C$ does not depend on $x_0, L$.
\end{proposition}
\proof
For any $-\f L2<a-1<a<b<b+1<\f L2$, define $\eta\in C^\infty (\ol)$, $0\leq\eta\leq 1$, $|\g\eta|\leq 2$ by
\[
\eta(x) =
\begin{cases}
0,\qquad x_n\leq a-1,\\
1,\qquad a\leq x_n\leq b,\\
0,\qquad x_n\geq b+1.
\end{cases}
\]
For any constants $k_1,k_2$, set
\[
\hat\fai(x) =
\begin{cases}
\fai(x) - k_1,&\qquad x_n\leq a,\\
\fai(x)-k_1 -\f{k_2-k_1}{b-a}(x_n-a),&\qquad a\leq x_n\leq b,\\
\fai(x) - k_2,&\qquad x_n\geq b.
\end{cases}
\]
Then $\eta^2\hat\fai\in H^1(\ol)$ and
$\left(\eta^2\hat\fai\right)|_{x_n=\pm L} = 0$. Therefore
$\eta^2\hat\fai\in H_L$ and
\[
\int_\ol \H(|\g\fai|^2)\g\fai\cdot\g(\eta^2\hat\fai)dx = 0.
\]
Thus,
\[
\int_{\O_{a-1,b+1}}\eta^2\H(|\g\fai|^2)\g\fai\cdot\g\hat\fai dx =
-2\int_{\O_{a-1,b+1}}\eta\H(|\g\fai|^2)\g\fai\cdot\g\eta\hat\fai dx,
\]
where $\g\hat\fai = \g\fai -
\f{k_2-k_1}{b-a}\chi_{a,b}(x)\vec{e}_n$, $\vec{e}_n=(0,...,0,1)$,
$$\O_{a,b}=\{x=(x_1,x_2,\cdots,x_n)\in\O|a<x_n<b\}$$ and
$\chi_{a,b}(x)$ is the characteristic function of $\O_{a,b}$. Then,
\[
\ba &\int_{\O_{a-1,b+1}}\eta^2\H(|\g\fai|^2)|\g\fai|^2 dx
+\int_{\O_{a,b}}\eta^2\H(|\g\fai|^2)\f{\p\fai}{\p x_n}\left(-\f{k_2-k_1}{b-a}\right)dx\\
=& -2\int_{\O_{a-1.b+1}}\eta\H(|\g\fai|^2)\g\fai\cdot\g\eta\hat\fai
dx. \ea
\]
Since $\eta =1$ on $\O_{a,b}$ and
$\disp\int_{S_{x_n}}\H(|\g\fai|^2)\f{\p\fai}{\p x_n}dx' = m_0$,
 \[
\int_{\O_{a-1,b+1}}\eta^2\H(|\g\fai|^2)|\g\fai|^2 dx =
-2\int_{\O_{a-1.b+1}}\eta\H(|\g\fai|^2)\g\fai\cdot\g\eta\hat\fai dx
+(k_2-k_1)m_0.
\]
Consequently, \be\label{n23} \ba
\ld\int_{\O_{a,b}}|\g\fai|^2dx &\leq \int_{\O_{a-1,b+1}}\eta^2\H(|\g\fai|^2)|\g\fai|^2 dx\\
&\leq 2\bigg|\int_{\O_{a-1,a}} +\int_{\O_{b,b+1}}\eta\H(|\g\fai|^2)\g\fai\cdot\g\eta\hat\fai dx\bigg| +|k_2-k_1|m_0\\
&\leq 4\Ld\bigg|\int_{\O_{a-1,a}} +\int_{\O_{b,b+1}}|\g\fai||\hat\fai|dx\bigg| +|k_2-k_1|m_0\\
&\leq 4\Ld\bigg[\bigg(\int_{\O_{a-1,a}}|\g\fai|^2dx\bigg)^{\f 12}\bigg(\int_{\O_{a-1,a}}|\fai - k_1|^2dx\bigg)^{\f 12}\\
&\ \ + \bigg(\int_{\O_{b,b+1}}|\g\fai|^2dx\bigg)^{\f
12}\bigg(\int_{\O_{b,b+1}}|\fai - k_2|^2dx\bigg)^{\f 12}\bigg]
+|k_2-k_1|m_0. \ea \ee Set
\[
k_1 = \aint_{\O_{a-1,a}}\fai dx,\qquad k_2 = \aint_{\O_{b,b+1}}\fai
dx.
\]
It follows the uniform Poincar\'{e} Inequality that \be\label{n24}
\int_{\O_{a-1,a}}|\fai - k_1|^2dx\leq
C\int_{\O_{a-1,a}}|\g\fai|^2dx,\ \ \int_{\O_{b,b+1}}|\fai -
k_2|^2dx\leq C\int_{\O_{b,b+1}}|\g\fai|^2dx, \ee where $C$ does not
depend on $a,b$.

\noindent Therefore, substituting (\ref{n24}) into (\ref{n23}) yields
\[
\ld\int_{\O_{a,b}}|\g\fai|^2dx \leq C\Ld\int_{\O_{a-1,b+1}\backslash\O_{a,b}}|\g\fai|^2dx + |k_2-k_1|m_0.
\]

We now claim that
\be\label{eqclaim}
|k_2-k_1|\leq C\int_{\O_{a-1,b+1}}|\g\fai|dx.
\ee

\noindent Assuming (\ref{eqclaim}) for a moment, one gets
\be\label{n25} \int_{\O_{a,b}}|\g\fai|^2 dx \leq
C\f{\Ld}{\ld}\int_{\O_{a-1,b+1}\backslash\O_{a,b}}|\g\fai|^2dx + \f
C\ld m_0\int_{\O_{a-1,b+1}}|\g\fai|dx.
\ee
Set $S_{\max} = \max_{x_n} |S_{x_n}|$, $S_{\min}=\min_{x_n} |S_{x_n}|$. On another hand,
\be\label{n26} \ba m_0\int_{\O_{a-1,b+1}}|\g\fai|dx &\leq
m_0S_{max}^{\f12}\bigg(\int_{\O_{a-1,b+1}}|\g\fai|^2dx\bigg)^{\f 12}
(b-a+2)^{\f 12}\\
&\leq \eps \int_{\O_{a-1,b+1}}|\g\fai|^2dx + \f
{S_{max}}{4\eps}(b-a+2)m_0^2. \ea\ee Combining (\ref{n25}) and
(\ref{n26}) leads to
\[\ba
\int_{\O_{a,b}}|\g\fai|^2 dx &\leq
C\f{\Ld}{\ld}\int_{\O_{a-1,b+1}\backslash\O_{a,b}}|\g\fai|^2dx + \f
C\ld \eps \int_{\O_{a-1,b+1}}|\g\fai|^2dx\\ &\ \ \ \ + \f
{C'}{\eps\ld}(b-a+2)m_0^2.\ea
\]
Taking $\disp \f{C\eps}\ld = \f 12$ yields
\[
\f 12\int_{\O_{a,b}}|\g\fai|^2 dx \leq \left(C\f{\Ld}{\ld} +\f
12\right)\int_{\O_{a-1,b+1}\backslash\O_{a,b}}|\g\fai|^2dx +  \f
{C'}{\ld^2}(b-a+2)m_0^2.
\]
Therefore, one has
\[
\int_{\O_{a,b}}|\g\fai|^2 dx \leq \left(C\f{\Ld}{\ld}
+1\right)\int_{\O_{a-1,b+1}\backslash\O_{a,b}}|\g\fai|^2dx +  \f
{C'}{\ld^2}(b-a+2)m_0^2,
\]
ie.
\[
\left(C\f{\Ld}{\ld} +2\right)\int_{\O_{a,b}}|\g\fai|^2 dx \leq
\left(C\f{\Ld}{\ld} +1\right)\int_{\O_{a-1,b+1}}|\g\fai|^2dx +  \f
{C'}{\ld^2}(b-a+2)m_0^2.
\]
Set $\disp \th_0 = \f{C\f{\Ld}{\ld} +1}{C\f{\Ld}{\ld} +2}$. Then
$0<\th_0<1$ and \be\label{ii1} \int_{\O_{a,b}}|\g\fai|^2 dx \leq
\th_0\int_{\O_{a-1,b+1}}|\g\fai|^2dx +  \f{C'}{\ld(C\Ld
+2\ld)}(b-a+2)m_0^2. \ee Set
\[
A_{a,b} = \f 1{b-a}\int_{\O_{a,b}}|\g\fai|^2dx.
\]
It follows from (\ref{ii1}) that
\[ A_{a,b}
\leq\th_0\f{b-a+2}{b-a}A_{a-1,b+1} + \f {C'}{\ld(C\Ld
+2\ld)}\f{b-a+2}{b-a}m_0^2.
\]
Taking $\theta'_0=\frac{1+\theta_0}2$ and a positive constant
$k(\th_0)\geq 2$ such that, if $b-a\geq k(\th_0)$, one has
\[
\disp\f{b-a+2}{b-a} \leq2,\qquad \disp\th_0\f{b-a+2}{b-a}\leq
\th_0'<1,
\]
\[
A_{a,b}\leq \th_0'A_{a-1,b+1} + \f {C'}{\ld(C\Ld +2\ld)}m_0^2\quad\
\ \text{for}\ \ \forall b-a\geq k(\th_0).
\]
Then,
\[\begin{array}{rl}
A_{a,b}&\leq (\th_0')^N A_{a-N,b+N} + \f {Cm_0^2}{\ld(C\Ld
+2\ld)}\sum_{i=0}^{N-1}(\th_0')^i\\
&\leq (\th_0')^N A_{a-N,b+N} + \f
{C'm_0^2}{\ld(C\Ld +2\ld)(1-\th_0')}.
\end{array}
\]
Applying (\ref{L}), one has
\[\begin{array}{rl}
A_{a,b}&\leq (\th_0')^N \f{|\O_{a-N,b+N}|}{b-a+2N}Cm_0^2 + \f
{Cm_0^2}{\ld(C\Ld +2\ld)(1-\th_0')}\\
&\leq C(\th_0')^N S_{max}m_0^2 +
\f {C'm_0^2}{\ld(C\Ld +2\ld)(1-\th_0')}.\end{array}
\]

Therefore, for any $\disp -\f 23L<a<b<\f 23 L$ and $b-a\geq
k(\th_0)$, letting $N \rightarrow \infty$ yields
\[
A_{a,b}\leq C m_0^2,
\]
where $C$ does not depend on $L$.

Then, for any $x_0$, $\disp |x_{0,n}|\leq \f L2$,
\[
\ba
\f 1{|\O_{x_0,1}|}\int_{\O_{x_0,1}}|\g\fai|^2 dx &\leq \f 1{|\O_{x_0,1}|}\int_{\O_{x_0,k(\th_0)}}|\g\fai|^2dx\\
&\leq  \f {2k(\th_0)S_{\max}}{|\O_{x_0,1}|}\cdot\f
1{2k(\th_0)S_{\max}}\int_{\O_{x_0,k(\th_0)}}|\g\fai|^2dx\\
&\leq \f {2k(\th_0)S_{\max}}{|\O_{x_0,1}|} Cm_0^2\\
&\leq  C\f {2k(\th_0)S_{\max}}{S_{\min}}m_0^2, \ea
\]
which yields (\ref{eq-LAE}).

Now, it remains to prove the claim (\ref{eqclaim}). For any $a\in
\R$, we define
$$
\a_i = \aint_{\O_{a+i-1,a+i}}\fai dx,\ \ \ \ \a_{i+\f 12} = \aint_{\O_{a+i-1,a+i+1}}\fai dx,\ \ \ i=1,2,\cdots, n.
$$
By the uniform Poincar\'{e} inequality (\ref{eq-upi}), one has
\be\label{n27}\ba \int_{\O_{a-1,a}}|\fai-\a_0|dx\leq
C\int_{\O_{a-1,a}}|\g\fai|dx\leq C\int_{\O_{a-1,a+1}}|\g\fai|dx, \ea
\ee
\be\label{n28}
\ba\int_{\O_{a-1,a}}\left|\fai-\a_{\f
12}\right|dx\leq\int_{\O_{a-1,a+1}}\left|\fai-\a_{\f
12}\right|dx\leq C\int_{\O_{a-1,a+1}}|\g\fai|dx. \ea\ee Then, it
follows (\ref{n27}) and (\ref{n28}) that
\[
\int_{\O_{a-1,a}}\left|\a_{\f 12}-\a_0\right|dx\leq
C\int_{\O_{a-1,a+1}}|\g\fai|dx.
\]
Consequently,
\[
\left|\a_{\f 12}-\a_0\right|\leq \f
C{|\O_{a-1,a}|}\int_{\O_{a-1,a+1}}|\g\fai|dx\leq \f
C{S_{\min}}\int_{\O_{a-1,a+1}}|\g\fai|dx.
\]
Similarly,
\[
\left|\a_1-\a_{\f 12}\right|\leq \f
C{S_{\min}}\int_{\O_{a-1,a+1}}|\g\fai|dx.
\]
Hence,
\[
|\a_1-\a_0|\leq \left|\a_1-\a_{\f 12}\right| + \left|\a_{\f
12}-\a_0\right|\leq \f C{S_{\min}}\int_{\O_{a-1,a+1}}|\g\fai|dx.
\]
In a similar way, one gets
\[
|\a_2-\a_1|\leq \f C{S_{\min}}\int_{\O_{a,a+2}}|\g\fai|dx,\ \ \ldots
,\ \ |\a_{n}-\a_{n-1}|\leq \f
C{S_{\min}}\int_{\O_{a+n-2,a+n}}|\g\fai|dx.
\]
Therefore, by induction, it holds that
\[
|\a_n - \a_0|\leq \f C{S_{\min}}\int_{\O_{a-1,a+n}}|\g\fai|dx,
\]
which proves the claim (\ref{eqclaim}).
\endproof

\subsection{$C^{1,\a}$ regularity of the weak solution}
\begin{lemma}\label{lem-ge} %%gradient estimate
(Gradient estimate). It holds that \be\label{eq-ge}
\|\g\fai\|_{L^{\infty}\left(\O_{L/2}\right)}\leq C m_0. \ee where
$C$ does not depend on $L$.
\end{lemma}
\proof \noindent The proof is based on Moser's iteration technique.\\
{\bf Step 1. Interior
estimate:}\def\ro{2R}\def\ri{R}\def\bo{{B_{\ro}}}\def\bi{B_{\ri}} It
follows from the definition of weak solutions that for any
$\bo\subset \ol$
\[
\int_\bo \H(|\g\fai|^2)\g\fai\cdot\g\psi dx = 0,\quad\forall \psi\in C^\infty_0(\bo).
\]
Regarding $\p_s\psi$ as a test function, $s=1,2,\cdots,n$, one gets
\[\ba
0 &= \int_\bo \H(|\g\fai|^2)\g\fai\cdot\g(\p_s\psi) dx
 = -\int_\bo \p_s(\H(|\g\fai|^2)\g\fai)\cdot\g\psi dx\\
&= -\int_\bo \bigg(\H(|\g\fai|^2)\d_{ij} +2\H'(|\g\fai|^2)\p_i\fai\p_j\fai\bigg)\p_i(\p_s\fai)\p_j\psi dx\\
&= -\int_\bo a_{ij}\p_i w_s\p_j\psi dx, \ea\] where
\[
a_{ij} = \H(|\g\fai|^2)\d_{ij} +2\H'(|\g\fai|^2)\p_i\fai\p_j\fai \in
L^\infty(\bo),\quad w_s = \p_s\fai\in L^2(\bo).
\]
Therefore \be\label{eqwk1} \int_\bo a_{ij}\p_i w_s\p_j\psi dx =
0,\qquad \forall \psi\in H^1_0(\bo). \ee Taking
\[
\psi = \eta^2 w_s^{p-1},\quad \eta\in C^\infty_0(\bo),\ \ \
\eta\equiv 1 \ \ \text{in}\ B_R,\ \ \ p\geq 2.
\]
in (\ref{eqwk1}), one has
\[\ba
0&= \int_\bo a_{ij}\p_i w_s\p_j\left(\eta^2 w_s^{p-1}\right) dx\\
 &= (p-1)\int_\bo \eta^2 a_{ij}\p_i w_s\p_j w_s  w_s^{p-2} dx + 2\int_\bo \eta a_{ij}\p_i w_s\p_j \eta w_s^{p-1} dx.
\ea\] Therefore
\[\ba
(p-1)\int_\bo& \eta^2 a_{ij}\p_i w_s\p_j w_s  w_s^{p-2} dx
\leq 2\bigg|\int_\bo \eta a_{ij}\p_i w_s\p_j \eta w_s^{p-1} dx\bigg|\\
&\leq 2\bigg(\int_\bo \eta^2 a_{ij}\p_i w_s\p_j w_s w_s^{p-2} dx\bigg)^{\f 12} \bigg(\int_\bo a_{ij}\p_i \eta\p_j \eta w_s^{p} dx\bigg)^{\f 12}\\
&\leq \f{p-1}{2} \int_\bo \eta^2 a_{ij}\p_i w_s\p_j w_s w_s^{p-2} dx
+ \f{2}{p-1}\int_\bo a_{ij}\p_i \eta\p_j \eta w_s^{p} dx, \ea\] ie.
\[
\int_\bo \eta^2 a_{ij}\p_i w_s\p_j w_s  w_s^{p-2} dx\leq
\f{4}{(p-1)^2}\int_\bo a_{ij}\p_i \eta\p_j \eta w_s^{p} dx.
\]
Due to (\ref{fact2}), we have
\be\label{n29} \int_\bo \eta^2 |\g w_s|^2  w_s^{p-2} dx\leq
\f{4\Ld}{(p-1)^2\ld}\int_\bo |\g \eta|^2 w_s^{p} dx.
\ee
Since
\be\label{n30} \ba \eta^2|\g w_s|^2  w_s^{p-2}
&=\f 4{p^2}\left|\g\left(\eta w_s^{\f p2}\right) - w_s^{\f p2}\g
\eta\right|^2\geq \f 2{p^2}\left|\g\left(\eta w_s^{\f
p2}\right)\right|^2 - \f 4{p^2}|\g\eta|^2 w_s^p,
\ea\ee
 Combining (\ref{n29}) with (\ref{n30}) yields that
\[
\f 2{p^2}\int_\bo \left|\g\left(\eta w_s^{\f p2}\right)\right|^2
dx\leq \f 4{p^2}\int_\bo \left|\g\eta\right|^2 w_s^p dx +
\f{4\Ld}{(p-1)^2\ld}\int_\bo |\g \eta|^2 w_s^{p} dx,
\]namely,
\[
\int_\bo \left|\g\left(\eta w_s^{\f p2}\right)\right|^2 dx \leq
2\bigg(\f{p^2\Ld}{(p-1)^2\ld} +1\bigg)\int_\bo |\g \eta|^2 w_s^{p}
dx.
\]
Then the Sobolev's inequality implies that
\be\label{eqm1}
\ba
\bigg(\int_\bo \left(\eta w_s^{\f p2}\right)^{\f{2n}{n-2}}
dx\bigg)^{\f{n-2}{n}}
&\leq C\int_\bo \left|\g\left(\eta w_s^{\f
p2}\right)\right|^2 dx\\
&\leq C\bigg(\f{p^2\Ld}{(p-1)^2\ld} +1\bigg)\int_\bo |\g \eta|^2
w_s^{p} dx.
\ea\ee
 for $s=1,2,\cdots,n$.

Set
\[
\ba
&p_k = p\bigg(\f{n}{n-2}\bigg)^k,\quad R_k = \ri\left(1+\f 1{2^k}\right),\\
&\eta_k\in C^{\infty}_0(B_{R_k}),\ \ \eta_k\equiv1\text{ in
}B_{R_{k+1}},\ \ \ |\g\eta_k|^2\leq \f C{R_{k} - R_{k+1}} =
\f{C2^{k+1}}{\ri}. \ea
\]
Note that, $\{p_k\}$ is a strictly increasing sequence and tends to
infinity as $k\rightarrow+\infty$, and $\{R_k\}$ is strictly
decreasing sequence and tends to $\ri$ as $k$ goes to infinity. The
following is the standard Moser's iteration process.

Taking $p= p_k$, $\eta=\eta_k$ in (\ref{eqm1}) yields that
\[\ba
\bigg(\int_{B_{R_{k+1}}} w_s^{p_{k+1}} dx\bigg)^{\f{n-2}{n}} &\leq
C\bigg(\f{p_k^2\Ld}{(p_k-1)^2\ld} +1\bigg)\int_{B_{R_k}} |\g
\eta_k|^2 w_s^{p_k} dx\\
 &\leq C\f{\Ld 2^k}{\ld \ri}\int_{B_{R_k}}
w_s^{p_k} dx,\ea
\]
and so,
\[
\|w_s\|_{L^{p_{k+1}}\left(B_{R_{k+1}}\right)}\leq \bigg[C\f{\Ld
2^k}{\ld \ri}\bigg]^{\f 1{p_k}}
\|w_s\|_{L^{p_{k}}\left(B_{R_{k}}\right)}.
\]
Let $\disp M_k = \|w_s\|_{L^{p_{k}}\left(B_{R_{k}}\right)}$ and
$\disp D_k = \bigg[C\f{\Ld 2^k}{\ld \ri}\bigg]^{\f 1{p_k}}$. Then by
induction \be\label{n31} M_{k+1}\leq D_k M_k\leq\cdots\leq D_k
D_{k-1}\cdots D_0 M_0= M_0\prod_{j=0}^{k} D_j. \ee Due to
\[
\prod_{j=0}^{\infty} D_j = \prod_{j=1}^{\infty}\bigg[C\f{\Ld
2^j}{\ld \ri}\bigg]^{\f 1{p_j}} =\bigg[C\f{\Ld}{\ld
\ri}\bigg]^{\sum_{j=0}^{\infty}\f 1{p_j}} 2^{\sum_{j=0}^\infty
\f{j}{p_j}} =C\bigg(C\f{\Ld}{\ld \ri}\bigg)^{\f n{2p}},
\]
so taking $k\rightarrow\infty$ in (\ref{n31}) yields
\be\label{eqge1} \sup_{\bi} w_s \leq C\bigg(C\f{\Ld}{\ld
\ri}\bigg)^{\f n{2p}} \|w_s\|_{L^p(\bo)},\quad \text{for\ any}\
s=1,2,\cdots,n. \ee

\noindent{\bf Step 2. Boundary estimate:} For any $x_0\in
\p\O_{L/2}$, according to Lemma \ref{lem-bd}, there exists a
neighbourhood $U_{x_0}$ of $x_0$ in $\R^n$ and an invertible
$C^{2,\a}$ map $$T_{x_0}: U_{x_0}\cap \O_{L/2}\rightarrow
B_{\d_0}^+: x\mapsto y$$ satisfying (\ref{eq-bd}), where $B_{\d_0}$
is independent of $x_0$.
Define
\def\bd{{B^+_{\d_0}}}\def\bbd{{\p B^+_{\d_0}\cap\{y_n=0\}}}
\def\br{{B^+_R}}\def\bbr{{\p B^+_{R}\cap\{y_n=0\}}}
\[
\sigma_{ij} = \f{\p y_j}{\p x_i}.
\]
Then (\ref{eq-bd}) implies
\[
\sigma_{ij}(y)\sigma_{in}(y) = 0, \quad \text{ on }\ y_n = 0,\quad
j=1,2,\cdots,n-1.
\]
For any $0<R\leq\d_0$, \be\label{eqw2} \int_\br \H
\sigma_{ij}\t\p_j\t\fai \sigma_{il}\t\p_l\t\psi J dy=0,\qquad
\forall \t\psi\in H_0^1\left(\bd\right). \ee Taking ($s=1,2,\cdots,
n$)
\[
\t\psi =
\t\p_s\left(\t\eta^2\left(\t\p_s\t\fai\right)^{p-1}\right),\ \ \ \
\t\eta\in C_0^\infty\left(\bd\right)\ \ \ \text{and}\ \ \t\eta\equiv
1\ \ \text{in}\ B^{+}_{\frac{\delta_0}2}
\]
in (\ref{eqw2}) and integrating by parts show that
\begin{equation}\label{340}
\begin{array}{rl} 0
&= \int_\br \H \sigma_{ij}\t\p_j\t\fai \sigma_{il}\t\p_l\left(\t\p_s(\t\eta^2(\t\p_s\t\fai)^{p-1})\right)Jdy\\
\ &= -\int_\br \t\p_s\left(\H \sigma_{ij}\t\p_j\t\fai
\sigma_{il}J\right)\t\p_l\left(\t\eta^2(\t\p_s\t\fai)^{p-1}\right)dy\\
&~~
-\int_\bbr \left(\H \sigma_{ij}\t\p_j\t\fai \sigma_{il}\right)\t\p_l\left(\t\eta^2(\t\p_s\t\fai)^{p-1}\right)J \d_{sn}dy'\\
\ &= -\int_\br \t\p_s\left(\H \sigma_{ij}\t\p_j\t\fai
\sigma_{il}J\right)\t\p_l\left(\t\eta^2(\t\p_s\t\fai)^{p-1}\right)dy,
\end{array}
\end{equation}
where the boundary terms vanish according to (\ref{eq-bd}-2) and
(\ref{eq-wb}).

Detailed calculations show that
\[\ba
\t\p_s(\H \sigma_{ij}\t\p_j\t\fai \sigma_{il}J)
&= \H \sigma_{ij}\t\p_j(\t\p_s\t\fai) \sigma_{il}J + 2\H'\sigma_{\a\b}\t\p_\b\t\fai\t\p_s(\sigma_{\a\gamma}\t\p_{\gamma}\t\fai) \sigma_{ij}\t\p_j\t\fai \sigma_{il}J\\
&\ \ \ + \H \t\p_j\t\fai\t\p_s(\sigma_{ij} \sigma_{il}J) \\
&= \H \sigma_{ij}\t\p_j(\t\p_s\t\fai) \sigma_{il}J
 + 2\H'\sigma_{\a\b}\t\p_\b\t\fai \sigma_{\a\gamma}\t\p_{\gamma}(\t\p_s\t\fai) \sigma_{ij}\t\p_j\t\fai \sigma_{il}J\\
&\qquad + 2\H'\sigma_{\a\b}\t\p_\b\t\fai\t\p_s(\sigma_{\a\gamma})\t\p_{\gamma}\t\fai \sigma_{ij}\t\p_j\t\fai \sigma_{il}J
 + \H \t\p_j\t\fai\t\p_s(\sigma_{ij} \sigma_{il}J)\\
&= \big(\H \sigma_{ij}\t\p_j(\t\p_s\t\fai)
 + 2\H'\sigma_{\a\b}\t\p_\b\t\fai \sigma_{\a\gamma}\t\p_{\gamma}(\t\p_s\t\fai) \sigma_{ij}\t\p_j\t\fai\big) \sigma_{il}J\\
&\qquad + 2\H'\sigma_{\a\b}\t\p_\b\t\fai\t\p_s(\sigma_{\a\gamma})\t\p_{\gamma}\t\fai \sigma_{ij}\t\p_j\t\fai \sigma_{il}J
 + \H \t\p_j\t\fai\t\p_s(\sigma_{ij} \sigma_{il}J)\\
&= \big(\H \d_{i\a}
 + 2\H'\sigma_{\a\b}\t\p_\b\t\fai  \sigma_{ij}\t\p_j\t\fai\big) \sigma_{\a\gamma}\t\p_{\gamma}(\t\p_s\t\fai) \sigma_{il}J\\
&\qquad + 2\H'\sigma_{\a\b}\t\p_\b\t\fai\t\p_s(\sigma_{\a\gamma})\t\p_{\gamma}\t\fai \sigma_{ij}\t\p_j\t\fai \sigma_{il}J
 + \H \t\p_j\t\fai\t\p_s(\sigma_{ij} \sigma_{il}J)\\
&= \big(\H \d_{ij}
 + 2\H'\sigma_{i\a}\t\p_\a\t\fai  \sigma_{j\b}\t\p_\b\t\fai\big) \sigma_{i\gamma}\t\p_{\gamma}(\t\p_s\t\fai) \sigma_{jl}J\\
&\qquad + 2\H'\sigma_{\a\b}\t\p_\b\t\fai\t\p_s(\sigma_{\a\gamma})\t\p_{\gamma}\t\fai \sigma_{ij}\t\p_j\t\fai \sigma_{il}J
 + \H \t\p_j\t\fai\t\p_s(\sigma_{ij} \sigma_{il}J).
\ea
\]
Set
\[
\ba
\t A_{l\gamma} &= \big(\H \d_{ij}
 + 2\H'\sigma_{i\a}\t\p_\a\t\fai  \sigma_{j\b}\t\p_\b\t\fai\big) \sigma_{i\gamma} \sigma_{jl}J, \\
\t B_{ls} & = 2\H'\sigma_{\a\b}\t\p_\b\t\fai\t\p_s(\sigma_{\a\gamma})\t\p_{\gamma}\t\fai \sigma_{ij}\t\p_j\t\fai \sigma_{il}J
 + \H \t\p_j\t\fai\t\p_s(\sigma_{ij} \sigma_{il}J).
\ea
\]
Then
\be\label{1eqwkb}
\int_\br \t
A_{l\gamma}\t\p_\gamma(\t\p_s\t\fai)\t\p_l(\t\eta^2(\t\p_s\t\fai)^{p-1})
+ \t B_{ls} \t\p_l(\t\eta^2(\t\p_s\t\fai)^{p-1}) dy =0.
\ee
Denoting $\t\p_s\t\fai$  by $\t w_s$, we have \be
\int_\br \t A_{l\gamma}\t\p_\gamma \t w_s\t\p_l(\t\eta^2\t
w_s^{p-1}) + \t B_{ls} \t\p_l(\t\eta^2\t w_s^{p-1}) dy =0, \ee which can be rewritten as,
\be\label{n35} \ba 0&=\int_\br \bigg(\t A_{l\gamma}\t\p_\gamma \t
w_s + \t B_{ls}\bigg)
\bigg((p-1)\t\eta^2 \t w_s^{p-2}\t\p_l\t w_s + 2\t\eta\t\p_l\t\eta \t w_s^{p-1}\bigg) dy\\
&= (p-1)\int_\br \t A_{l\gamma}\t\p_\gamma \t w_s  \t\eta^2 \t
w_s^{p-2}\t\p_l\t w_s dy
+ 2\int_\br \t A_{l\gamma}\t\p_\gamma \t w_s  \t\eta\t\p_l\t\eta \t w_s^{p-1} dy\\
&\ \ \ +(p-1)\int_\br \t B_{ls} \t\eta^2 \t w_s^{p-2}\t\p_l\t w_s dy
+2\int_\br \t B_{ls}\t\eta\t\p_l\t\eta \t w_s^{p-1} dy\\
&= I_1 + I_2 + I_3 + I_4. \ea \ee Note that,
\be\label{n32} \ba
|I_2| &\leq 2\bigg|\int_\br\t\eta \t A_{l\gamma}\t\p_l\t\eta\t\p_\gamma \t w_s \t w_s^{p-1}dy\bigg|\\
&\leq \f{p-1}{2}\int_\br\t\eta^2\t A_{l\gamma}\t\p_l \t
w_s\t\p_\gamma \t w_s \t w_s^{p-2}dy + \f{2}{p-1}\int_\br \t
A_{l\gamma}\t\p_l\t\eta\t\p_\gamma\t\eta \t w_s^p dy,\ea \ee
\be\label{n33} \ba
|I_3| &\leq (p-1)\bigg|\int_\br\t\eta^2 \t B_{ls}\t\p_l \t w_s \t w_s^{p-2}dy\bigg|\\
&\leq (p-1)\f\ld 4\int_\br\t\eta^2|\t\g \t w_s|^2 \t w_s^{p-2} dy + \f{p-1}\ld \int_\br \t\eta^2 |\t B_{ls}|^2 \t w_s^{p-2} dy,\\
\ea\ee and \be\label{n34} \ba |I_4| &\leq 2\bigg|\int_\br\t\eta \t
B_{ls}\t\p_l\t\eta \t w_s^{p-1}dy\bigg| \leq
2\int_\br\t\eta|\t\g\t\eta||\t B_{ls}| \t w_s^{p-1}dy. \ea \ee
Therefore, substituting (\ref{n32}), (\ref{n33}) and (\ref{n34}) into
(\ref{n35}) yields that
\[
\ba \f{I_1}2&=\f{(p-1)}2\int_\br\t\eta^2 \t A_{l\gamma}\t\p_l \t
w_s\t\p_\gamma \t w_s
\t w_s^{p-2}dy\\
&\leq  \f{2}{p-1}\int_\br \t A_{l\gamma}\t\p_l\t\eta\t\p_\gamma\t\eta \t w_s^p dy + (p-1)\f\ld 4\int_\br\t\eta^2|\t\g \t w_s|^2 \t w_s^{p-2} dy\\
&\ \ \  + \f{p-1}\ld \int_\br \t\eta^2 |\t B_{ls}|^2 \t w_s^{p-2} dy +2\int_\br\t\eta|\t\g\t\eta||\t B_{ls}| \t w_s^{p-1}dy.\\
\ea\] Then the uniform ellipticity yields
\[\ba
\int_\br\t\eta^2 |\t\g \t w_s|^2 \t w_s^{p-2}dy&\leq
\f{8\Ld}{\ld(p-1)^2}\int_\br |\t\g\t\eta |^2 \t w_s^p dy +
\f{4}{\ld^2} \int_\br \t\eta^2 |\t B_{ls}|^2 \t w_s^{p-2} dy\\
&\ \ \ \  +\f {8}{\ld(p-1)}\int_\br\t\eta|\t\g\t\eta||\t B_{ls}| \t
w_s^{p-1}dy. \ea
\]
Note that (\ref{fact1}) implies that
\[\ba
|\t B_{ls}|&\leq C(1+\Ld)|\t w|\leq C\Ld |\t w|,
 \ea
\]
where $\t w = (\t w_1, \t w_2,\cdots, \t w_n)= (\t\p_1\t\fai,
\t\p_2\t\fai,\cdots, \t\p_n\t\fai)$ and $\disp |\t w| = \sum_{s=1}^n
|\t w_s|$. Then,
\[\ba
\int_\br\t\eta^2 |\t\g \t w_s|^2 \t w_s^{p-2}dy &\leq
\f{8\Ld}{\ld(p-1)^2}\int_\br |\t\g\t\eta |^2 \t w_s^p dy +
\f{C\Ld^2}{\ld^2} \int_\br \t\eta^2 |\t w|^2 \t w_s^{p-2} dy\\
&\ \ \ \ +\f {C\Ld}{\ld(p-1)}\int_\br\t\eta|\t\g\t\eta||\t w| \t
w_s^{p-1}dy \ea
\]
which implies
\be\label{n36}\ba \int_\br\t\eta^2 \left|\t\g\left(\t
w_s^{\f p2}\right)\right|^2dy&\leq \f{2\Ld p^2}{\ld(p-1)^2}\int_\br
|\t\g\t\eta |^2 \t w_s^p dy + \f{C\Ld^2}{\ld^2}p^2 \int_\br \t\eta^2
|\t w|^2 \t w_s^{p-2} dy\\& \ \ \ \ +\f {C\Ld
p^2}{\ld(p-1)}\int_\br\t\eta|\t\g\t\eta||\t w| \t w_s^{p-1}dy. \ea
\ee It follows from (\ref{n36}) that
\be\label{eqm2} \ba &\ \ \ \int_\br \left|\t\g\left(\t\eta \t w_s^{\f
p2}\right)\right|^2dy\\
& \leq \bigg(\f{2\Ld
p^2}{\ld(p-1)^2}+1\bigg)\int_\br |\t\g\t\eta |^2 \t w_s^p dy
+ \f{C\Ld^2}{\ld^2}p^2 \int_\br \t\eta^2 |\t w|^2 \t w_s^{p-2} dy\\
&\ \ \ \ \ +\f {C\Ld p^2}{\ld(p-1)}\int_\br\t\eta|\t\g\t\eta||\t w|
\t w_s^{p-1}dy, \ea \ee for $s=1,2,\cdots,n$.

Let
\[
\ba &R_k = \d_0\left(\th +\f{1-\th}{2^k}\right),\quad \ \ \ p_k =
p\left(\f{n}{n-2}\right)^k,\ \ \ \ \ k=0,1,2,\ldots,\\ & \t\eta_k\in
C^{\infty}_0(B_{R_k}),\ \ \t\eta_k\equiv 1\text{ in }B_{R_{k+1}}\ \
\ \text{and}\ \
|\t\g\t\eta_k|\leq \f{2}{R_k-R_{k+1}}=\f{2^{k+1}}{(1-\th)\d_0}.\\
\ea\] Taking $p= p_k$, $\t\eta=\t\eta_k$ in (\ref{eqm2}) and using Sobolev embedding Theorem, we obtain
\[
\ba &\ \ \ \ \ \bigg(\int_{B^+_{R_{k+1}}} \t w_s^{p_{k+1}}
dy\bigg)^{\f{n-2}{n}}\\
&\leq \bigg(\int_{B^+_{R_{k}}}\left(\t\eta_k \t
w_s^{\f{p_{k}}{2}}\right)^{\f{2n}{n-2}}dy\bigg)^{\f{n-2}{n}}\\
&\leq C\int_{B^+_{R_{k}}} \left|\t\g\left(\t\eta_k \t w_s^{\f {p_k}2}\right)\right|^2dy\\
&\leq  C\bigg(\f{2\Ld
{p_k}^2}{\ld({p_k}-1)^2}+1\bigg)\int_{B^+_{R_{k}}} |\t\g\t\eta_k |^2
\t w_s^{p_k} dy+ \f{C\Ld^2}{\ld^2}{p_k}^2 \int_{B^+_{R_{k}}}
\t\eta_k^2
|\t w|^2 \t w_s^{{p_k}-2} dy\\
&\ \ \ \ +\f {C\Ld {p_k}^2}{\ld(p_k-1)}\int_{B^+_{R_{k}}}\t\eta_k|\t\g\t\eta_k||\t w| \t w_s^{{p_k}-1}dy\\
&\leq  C\bigg(\f{2\Ld
{p_k}^2}{\ld({p_k}-1)^2}+1\bigg)\int_{B^+_{R_{k}}}
\f{4^{k+1}}{(1-\th)^2\d_0^2} \t w_s^{p_k} dy +
\f{C\Ld^2}{\ld^2}{p_k}^2
\int_{B^+_{R_{k}}}  |\t w|^2 \t w_s^{{p_k}-2} dy\\
&\ \ \ \  +\f {C\Ld {p_k}^2}{\ld({p_k}-1)}\int_{B^+_{R_{k}}}
\f{2^{k+1}}{(1-\th)\d_0}|\t w| \t w_s^{{p_k}-1}dy. \ea
\]
Therefore,
$$ \|\t w_s\|_{L^{p_{k+1}}\left(B_{R_{k+1}}^+\right)}
\leq  \bigg(A_k\int_{B^+_{R_{k}}}  \t w_s^{p_k} dy + B_k\int_{B^+_{R_{k}}}  |\t w|^2 \t
w_s^{{p_k}-2} dy + C_k \int_{B^+_{R_{k}}}
|\t w| \t w_s^{{p_k}-1}dy\bigg)^{\f 1{p_k}},
$$
with
$$
A_k = C\bigg(\f{2\Ld {p_k}^2}{\ld({p_k}-1)^2}+1\bigg)
\f{4^{k+1}}{(1-\th)^2\d_0^2},\ \ \  B_k =
\f{C\Ld^2}{\ld^2}{p_k}^2,$$ and $$C_k = \f {C\Ld
{p_k}^2}{\ld({p_k}-1)}\cdot\f{2^{k+1}}{(1-\th)\d_0}.
$$
Note that
\[
\ba
\int_{B^+_{R_k}} \t w_s^{p_k-2} |\t w|^2 dy &\leq \bigg(\int_{B^+_{R_k}}\t w_s^{p_k}dy\bigg)^{\f{p_k-2}{p_k}}\bigg(\int_{B^+_{R_k}}|\t
w|^{p_k}dy\bigg)^{\f 2{p_k}},\\
\int_{B^+_{R_k}} \t w_s^{p_k-1} |\t w| dy &\leq
\bigg(\int_{B^+_{R_k}}\t w_s^{p_k}dy\bigg)^{\f{p_k-1}{p_k}}
\bigg(\int_{B^+_{R_k}}|\t w|^{p_k}dy\bigg)^{\f 1{p_k}}. \ea
\]
Then
\[
\ba \|\t w_s\|_{L^{p_{k+1}}\left(B_{R_{k+1}}^+\right)} &\leq \bigg[A_k\|\t w_s\|_{L^{p_{k}}\left(B_{R_{k}}^+\right)}^{p_k}
+B_k \|\t w_s\|_{L^{p_{k}}\left(B_{R_{k}}^+\right)}^{p_k-2}\||\t w|\|_{L^{p_{k}}\left(B_{R_{k}}^+\right)}^{2}\\
&\quad +C_k \|\t
w_s\|_{L^{p_{k}}\left(B_{R_{k}}^+\right)}^{p_k-1}\||\t
w|\|_{L^{p_{k}}\left(B_{R_{k}}^+\right)}\bigg]^{\f
1{p_k}}, \\
&\leq \|\t
w_s\|_{L^{p_{k}}(B_{R_{k}}^+)}^{\f{p_k-2}{p_k}}\bigg[A_k\|\t
w_s\|_{L^{p_{k}}(B_{R_{k}}^+)}^2
+B_k \||\t w|\|_{L^{p_{k}}(B_{R_{k}}^+)}^{2}\\
&\quad +C_k \|\t w_s\|_{L^{p_{k}}(B_{R_{k}}^+)}\||\t
w|\|_{L^{p_{k}}(B_{R_{k}}^+)}\bigg]^{\f 1{p_k}}. \ea
\]
Hence,
\[
\ba \|\t w_s\|_{L^{p_{k+1}}\left(B_{R_{k+1}}^+\right)} &\leq \|\t
w_s\|_{L^{p_{k}}\left(B_{R_{k}}^+\right)}^{\f{p_k-2}{p_k}}\bigg[A_k
+ B_k +C_k\bigg]^{\f 1{p_k}} \||\t
w|\|_{L^{p_{k}}\left(B_{R_{k}}^+\right)}^{\f 2{p_k}}, \ea
\]
which implies that
\[
\ba \sum_{s=1}^{n}\|\t w_s\|_{L^{p_{k+1}}\left(B_{R_{k+1}}^+\right)}
&\leq \bigg[A_k + B_k +C_k\bigg]^{\f 1{p_k}} \||\t w|\|_{L^{p_{k}}\left(B_{R_{k}}^+\right)}^{\f 2{p_k}} \sum_{s=1}^n\|\t
w_s\|_{L^{p_{k}}\left(B_{R_{k}}^+\right)}^{\f{p_k-2}{p_k}}\\
&\leq \bigg[A_k + B_k +C_k\bigg]^{\f 1{p_k}} \||\t
w|\|_{L^{p_{k}}\left(B_{R_{k}}^+\right)}^{\f 2{p_k}}
\bigg(\sum_{s=1}^n\|\t
w_s\|_{L^{p_{k}}\left(B_{R_{k}}^+\right)}^{\f{p_k-2}{p_k}\cdot\f{p_k}{p_k-2}}\bigg)^{\f{p_k-2}{p_k}}
n^{\f{2}{p_k}}\\
&\leq \bigg[(A_k + B_k +C_k)n^2\bigg]^{\f 1{p_k}}
\bigg(\sum_{s=1}^n\|\t
w_s\|_{L^{p_{k}}\left(B_{R_{k}}^+\right)}\bigg)^{\f{2}{p_k}}\\
&\ \ \ \times \bigg(\sum_{s=1}^n\|\t
w_s\|_{L^{p_{k}}\left(B_{R_{k}}^+\right)}\bigg)^{\f{p_k-2}{p_k}}.
\ea
\]
Therefore,
\[
\sum_{s=1}^{n}\|\t w_s\|_{L^{p_{k+1}}\left(B_{R_{k+1}}^+\right)}
\leq \bigg[(A_k + B_k +C_k)n^2\bigg]^{\f 1{p_k}} \sum_{s=1}^n\|\t
w_s\|_{L^{p_{k}}\left(B_{R_{k}}^+\right)}.
\]
Define
\[
M_k = \sum_{s=1}^n\|\t
w_s\|_{L^{p_{k}}\left(B_{R_{k}}^+\right)},\quad D_k = \bigg[(A_k +
B_k +C_k)n^2\bigg]^{\f 1{p_k}}.
\]
Then
\[
M_{k+1}\leq D_k M_k.
\]
It is clear that $\disp p_k = p\left(\f{n}{n-2}\right)^k\leq p4^k$,
for $n\geq 3$, so
\[
\ba A_k + B_k +C_k &\leq C\bigg(\f{2\Ld p_k^2}{\ld(p_k-1)^2}+1\bigg)
\f{4^{k+1}}{(1-\th)^2\d_0^2} +\f{C\Ld^2}{\ld^2}p_k^2
+\f {C\Ld p_k^2}{\ld(p_k-1)}\f{2^{k+1}}{(1-\th)\d_0}\\
&\leq C\bigg(\f{\Ld}{\ld}+1\bigg) \f{4^{k}}{(1-\th)^2\d_0^2}
+C\f{\Ld^2}{\ld^2}p^2 16^k
+C\f {\Ld}{\ld}p\f{8^{k}}{(1-\th)\d_0}\\
&\leq T\cdot 16^k, \ea
\]
where $\disp T = C\bigg[\f{1}{(1-\th)^2\d_0^2}\f{\Ld}{\ld}
+\f{\Ld^2}{\ld^2}p^2 +\f {\Ld}{\ld}\f{p}{(1-\th)\d_0}\bigg]$, $C$
does not depend on $k$.

\noindent Then,
\[
\ba
M_{k+1}&\leq D_k M_k \leq D_k\cdot D_{k-1}\cdots D_0\cdot M_0\\
&\leq T^{\sum_{i=0}^k\f{1}{p_i}}\cdot
16^{\sum_{i=0}^k\f{i}{p_i}}\cdot M_0. \ea
\]
Note that
\[
\sum_{i=0}^\infty\f{1}{p_i} = \f n{2p},\qquad  \ \text{and}\ \ \ \ \
\  \sum_{i=0}^k\f{i}{p_i}<\infty.
\]
One has
\[
M_{k+1}\leq CT^{\f n{2p}} M_0,\qquad \forall k>0.
\]
Letting $k\rightarrow\infty$ shows that \be\label{eqge2}
\sum_{s=1}^n\sup\limits_{B_{\th \d_0}} |\t w_s|\leq CT^{\f
n{2p}}\sum_{s=1}^n \|\t w_s\|_{L^p\left(B^+_{\d_0}\right)}. \ee
Combining the interior estimate (\ref{eqge1}) with the boundary
estimate (\ref{eqge2}) yields the desired gradient estimate (\ref{eq-ge}).
\endproof

\begin{remark}
It has been assumed that $w_s \geq 0$ and $w_s$ is bounded in the
above proof. The boundness
 assumption could be eliminated by a standard technique
(see chapter 8 of \cite{GT}). If $w_s$ is not positive, we can repeat the
proof for $w_s^+$ and $w_s^-$ respectively.
\end{remark}

\begin{remark}
In the case that $n =2$, choosing $p_1=\infty$, one can obtain the
estimate similarly to (\ref{eqge2}).
\end{remark}
\begin{lemma}\label{th-ghe} %%gradient holder estimate
(H\"older estimate of gradient.) $\fai\in
C^{1,\a}\left(\cl\O_{L/2}\right)$ and \be\label{eqghe}
\|\g\fai\|_{C^{0,\a}\left(\O_{L/2}\right)}\leq C m_0, \ee where $C$
does not depend on $L$.
\end{lemma}
\proof {\bf Step 1. Interior Estimate.} For any $B_{2R}\subset\O$,
$w_s = \p_s\fai$ $(s=1,2,\cdots,n)$ is a weak solution to
\[
\p_i(a_{ij}\p_j w_s) = 0,
\]
in the sense of (\ref{eqwk1}), where $\disp a_{ij} = \rho(|\g\fai|^2)\d_{ij}
+2\rho'(|\g\fai|^2)\p_i\fai\p_j\fai$. Then, the desired interior H\"{o}lder estimate for
$w_s$ is just the standard interior H\"older estimate for the weak
solutions to second order elliptic equation with bounded
coefficients.\\
{\bf Step 2. Boundary Estimate.} Similar to (\ref{1eqwkb}), one has
for any $s=1,2,\cdots,n-1$, \be\label{2eqwkb} \int_\br \t
A_{l\gamma}\t\p_\gamma \t w_s\t\p_l \t\psi + \t B_{ls} \t\p_l \t\psi dy
=0, \quad \t\psi\in H_0^1(\br). \ee where $\br$, $\t A_{l\gamma}$
and $\t B_{ls}$ are same as in (\ref{1eqwkb}).

By an even symmetrizing procedure, $\hat w_s$, $\hat A_{l\gamma}$
and $\hat B_{ls}$ denote the even extensions of $\t w_s$, $\t
A_{l\gamma}$ and $\t B_{ls}$, respectively. Then $\hat w_s$ satisfies
\be\label{2eqwk2} \int_{B_R} \hat A_{l\gamma}\t\p_\gamma \hat
w_s\t\p_l \hat\psi + \hat B_{ls} \t\p_l \hat\psi dy =0, \quad
\hat\psi\in H_0^1(B_R), \ 1\leq s\leq n-1. \ee Since $\|\t
w\|_{L^{\infty}}\leq Cm_0$,
\[
\|\t A_{l\gamma}, \t B_{ls}\|_{L^{\infty}}\leq C,\quad  \|\hat
A_{l\gamma}, \hat B_{ls}\|_{L^{\infty}}\leq C.
\]
Therefore, for $1\leq s\leq n-1$, the standard interior De Giorgi estimate gives
\[
\|\t w_s\|_{C^\a \left(B^+_{R/2}\right)}\leq C\left(\|\t
w_s\|_{L^2(\br)} + \f 1\ld \|\t B_{ls}\|_{L^q(\br)}\right) \leq
Cm_0,\quad q>n.
\]
Now, we estimate $\t w_n$. For any $y_0\in B^+_{R/2}$, $r\leq
\max\left\{\f 16 R, \f 12\right\}$, taking $\t\psi = \t\eta^2 (\t
w_s-\bar w_s)$ in (\ref{2eqwk2}), $\t\eta\in
C_0^{\infty}(B_{2r}(y_0))$, $\t\eta\equiv 1$ in $B_r(y_0)$,
$|\t\g\t\eta|\leq \f 2 r$ and $\bar w_s = \aint_{B_{2r}(y_0)}\t w_s
dy$, one has
\[
\int_\br \t A_{l\gamma}\t\p_\gamma \t w_s\t\p_l (\t\eta^2 (\t w_s-
\bar w_s)) + \t B_{ls} \t\p_l (\t\eta^2 (\t w_s-\bar w_s)) dy =0.
\]
Therefore,
\[
\ba &\ \ \ \ \int_{\br}\t A_{l\gamma} \t\p_l\t w_s\t\p_\gamma\t w_s\t\eta^2 dy\\
&\leq \left|\int_{\br}\t B_{ls} \t\p_l\t w_s\t\eta^2 dy\right| + \left|2\int_{\br}\t A_{l\gamma}\t\p_\gamma\t w_s\t\eta\t\p_l\t\eta(\t w_s-\bar w_s)
dy\right| +\left|2\int_{\br}\t B_{ls}\t\eta \t\p_l\t\eta(\t w_s-\bar w_s)dy\right| \\
&\leq \frac1{\ld}\int_{\br}|\t B_{ls}|^2\t\eta^2dy +\f\ld4
\int_{\br}|\t\g \t w_s|^2\t\eta^2dy + \f{\ld}{4\Ld}\int_{\br}\t A_{l\gamma}\t\p_l\t w_s\t\p_\gamma\t w_s\t\eta^2 dy\\
&\ \ \ \ +\f{4\Ld}{\ld}\int_{\br}\t
A_{l\gamma}\t\p_l\t\eta\t\p_\gamma\t\eta(\t w_s-\bar w_s)^2dy
+\int_{\br}|\t B_{ls}|^2\t\eta^2dy
+\int_{\br}\left|\t\g\t\eta\right|^2(\t w_s-\bar w_s)^2 dy. \ea
\]
As a consequence,
\[
\int_{\br}\t\eta^2|\t\g \t w_s|^2 dy \leq C\bigg(\int_{\br\cap
B_{2r}(y_0)} |\t\g\t\eta|^2|\t w_s-\bar w_s|^2 dy + \int_{\br\cap
B_{2r}(y_0)} \t\eta^2|\t B_{ls}|^2dy\bigg).
\]
Noting that $\|\t w_s\|_{C^\a\left(B^+_{R/2}\right)}\leq Cm_0$ and
$|\t B_{ls}|\leq Cm_0$, one has that, for any $1\leq s\leq n-1$,
\be\label{n37} \ba \int_{\br\cap B_{r}(y_0)}|\t\g \t w_s|^2 &\leq
C\bigg(\int_{\br\cap B_{2r}(y_0)} |\t\g\t\eta|^2|\t w_s-\bar w_s|^2
dy + \int_{\br\cap
B_{2r}(y_0)} \t\eta^2|\t B_{ls}|^2dy\bigg)\\
&\leq Cm_0^2\left(r^{n-2+2\a} + r^n\right)\\
&\leq Cm_0^2 r^{n-2+2\a}.
\ea \ee
 According to the equation of $\t\fai$, one has
\be\label{n38}\ba \int_{\br\cap B_{r}(y_0)}|\t D_{nn}^2\t\fai|^2dy
&\leq C\bigg(\sum_{k=1}^{n-1}\int_{\br\cap B_{r}(y_0)}|\t\g \t w_s|^2dy + \int_{\br\cap B_{r}(y_0)}|\t\g\t\fai|^2dy\bigg)\\
&\leq C m_0^2 (r^{n-2+2\a} + r^n)\\
&\leq Cm_0^2 r^{n-2+2\a}. \ea\ee
Then, due to (\ref{n37}) and (\ref{n38}), one has by Theorem
\ref{th-morrey} that
\[
\|w\|_{C^\a\left(B^+_{R/2}\right)}\leq Cm_0.
\]
Now, the H\"{o}lder estimate of $\g\fai$ follows from Step 1 and
Step 2.
\endproof

\subsection{Proof of the existence of subsonic flows}
{\bf Proof of the statement (i) of Theorem 1.} For any fixed suitably
large $L$, according to previous subsections, one can get a $H^1$
function $\fai_L(x)$ such that $(\fai_L(x) - \fai_L(0))\in H_L$ is a
weak solution  to problem \ref{ptnozzle}. Set $\hat\fai_L(x) =
\fai_L(x) - \fai_L(0)$. Moreover, $\hat\fai_L\in
C^{1,\a}(\cl\O_{L/2})$ and
\[
\|\g\hat\fai_L\|_{C^{0,\a}\left(\O_{L/2}\right)}\leq C m_0.
\]
For any fixed $K\gg 1$,
if $L>2K$,
\[
\|\hat\fai_L\|_{C^{1,\a}\left(\O_K\right)} \leq C,
\]
where $C$ does not depend on $L$, and $\hat\fai_L$ satisfies
\[
 \int_{\O_{K}} \H\left(|\g\hat\fai_L|^2\right)\g\hat\fai_L\cdot\g\psi dx = 0,\qquad \forall \psi\in C^\infty_0(\O_{K}).
\]
Since $\hat\fai_L\in H_L\cap C^{1,\a}(\O_K)$ satisfies the equation
(\ref{eq-W}), one can check easily that
\[
\int_{S_{x_{0}}}\H\left(|\g\hat\fai_L|^2\right)\f{\p\hat\fai_L}{\p x_n}
dx' = m_0,\quad\text{for any } x_0\in \O_{K}.
\]
By a standard diagonal argument, there exists a $\fai\in
C^{1,\a}(\O)$ and a subsequence $\hat\fai_{L_n}$ such that for any
$K$,
\[
\lim_{n\ra\infty}\|\hat\fai_{L_n} - \fai\|_{C^{1,\a}
\left(\O_K\right)} = 0.
\]
Therefore, one has
\[
 \int_{\O} \H(|\g\fai|^2)\g\fai\cdot\g\psi dx = 0,\qquad \forall \psi\in C^\infty_0(\O),
\]
and
\[
\int_{S_{x_0}}\H(|\g\fai|^2)\f{\p\fai}{\p x_n} dx' =
m_0,\quad\text{for any }x_0\in\O.
\]
It is clear that
\[
\fai\in C^{1,\a}(\O)\bigcap H^1_{loc}(\O)
\]
and
\[
\|\g\fai\|_{C^{0,\a}(\O)}\leq Cm_0.
\]
Similar to the previous subsections, one can prove that $\fai\in
H^2_{loc}(\O)$ and $\fai$ is a strong solution to \be\label{eq-a}
\left\{ \ba
&\left(\H(|\g\fai|^2)\d_{ij} + 2\H'(|\g\fai|^2)\p_i\fai\p_j\fai\right)\p^2_{ij}\fai = 0\qquad \text{ in }\O,\\
&\f{\p\fai}{\p \vec n} = 0,\qquad \text{ on }\p\O. \ea \right. \ee
By the standard regularity theory for second order elliptic
equations, one gets that $\fai\in C^{2,\a}_{loc}(\cl\O)$ is a
solution to (\ref{eq-a}) with the property
\[
\|\g\fai\|_{C^{1,\a}(\O)}\leq C m_0.
\]
Choose $m_0$ small enough such that $Cm_0\leq 1-2\hat\d_0$. Then the
subsonic truncation automatically disappears, so $\fai\in
C^{2,\a}_{loc}(\cl\O)$ is a smooth solution to the original Problem \ref{op}. This proves the first part of Theorem \ref{mainth}.
\endproof
\begin{remark}
In fact, we can conclude that $\fai\in C^{\infty}(\O)$ by the
standard bootstrap argument.
\end{remark}

\section{Uniqueness of the global subsonic flow}
\begin{theorem} (Uniqueness)
Suppose that $\Omega$ satisfies the assumptions (\ref{H1}), and $\varphi_k\ (k=1, 2)$ are uniformly subsonic solutions to
the following problem
\[\left
 \{\begin{array}{ll}
\text{div} \left(\rho(|\nabla \varphi_k|^2\right)\nabla \varphi_k)=0,\ \ \ \ \ \ & \text{in}\ \ \ \Omega,\\[3mm]
\frac{\partial\varphi_k}{\partial{\vec n}}=0,& \text{on}\ \ \ \p\O,
 \end{array}
 \right.
\]
associated with the same incoming mass flux $m_0$. Then
$$
\nabla \varphi_1 = \nabla \varphi_2,\ \ \ \ \ \ \ \  \ \ \ in\ \
\Omega.
$$
\end{theorem}
\proof Set $\varphi = \varphi_1 - \varphi_2$. Then $\varphi$
satisfies
\begin{equation}\label{n5}
\left
 \{\begin{array}{ll}
 \displaystyle{\partial_i(A_{ij}\partial_j\varphi)=0},\ \ \ \ \ \ & in\ \ \ \Omega,\\[3mm]
  \displaystyle{\frac{\partial\varphi}{\partial{\vec n}}=0,}& on\ \ \
  \p\O,
 \end{array}
 \right.
\end{equation}
where
$$
A_{ij}=\int_0^1\rho(\bar q^2)\delta_{ij}+2\rho'(\bar
q^2)(s\p_j\fai_1+(1-s)\p_j\fai_2)(s\partial_i\varphi_1+(1-s)\partial_i\varphi_2)ds,
$$
$$\displaystyle{\bar q^2 = |s\nabla\varphi_1 +
(1-s)\nabla\varphi_2|^2}.$$ Moreover, there exist two positive
constants $\lambda < \Lambda$, such that for any vector $\xi\in
\mathbb{R}^n$
\begin{equation}\label{n6}
\lambda |\xi|^2 < A_{ij} \xi_i\xi_j < \Lambda |\xi|^2.
\end{equation}
Let $\eta(x) = \eta(x_n)$ be a $C^\infty_0$ function satisfying
$$
\eta(x_n) \equiv1\ \  \text{for}\ \ |x_n|\leq L;\ \ \ \ \eta(x_n)
\equiv0\ \ \text{for}\ \ |x_n|\geq L+1, \ \ \ and\ \
|\eta'(x_n)|\leq 2,
$$
Denote $\displaystyle{\Omega_{a,b}=\{x=(x',x_n)\in\O|a\leq x_n \leq
b\}}$ and for $L > 0$
$$\hat\varphi(x)=
\left\{\begin{array}{ll}
 \displaystyle{\varphi(x)-\varphi_L^-},\ \ \ \ \ \ \ \ \ \ \ \ \ & x\in\Omega\cap\{x_n\leq-L\},\\[1mm]
  \displaystyle{\varphi(x)-\varphi_L^--\frac{\varphi_L^+-\varphi_L^-}{2L}(x_n+L)},\ \ \ & x\in\Omega\cap\{-L\leq x_n\leq L\},\\
  \displaystyle{\varphi(x)-\varphi_L^+}, \ &  x\in\Omega\cap\{x_n\geq L\},
 \end{array}
 \right.
$$
where $$\displaystyle{\varphi_L^- = \frac1{|\Omega_{ -L-1,
-L}|}\int_{\Omega_{- L-1, -L}}\varphi(x) dx}, \ \ \
\displaystyle{\varphi_L^+ = \frac1{|\Omega_{L, L+1}|}\int_{\Omega_{
L, L+1}}\varphi(x) dx}.
$$
Note that $\displaystyle{\nabla\hat\varphi=\nabla\varphi
-\frac{\varphi_L^+-\varphi^-_L}{2L}\chi_{-L,L}(x)\vec{e}_n}$,
$\vec{e}_n=(0, \cdots, 0, 1)$, $\chi_{-L,L}(x)$ is the
characteristic function of $\Omega_{-L,L}$.

Multiplying on the both sides of the first
equation in (\ref{n5}) by $\eta^2 \hat\varphi$, and integrating it over $\Omega$, one
obtains
$$
\int_{\Omega_{-L-1,L+1}}\eta^2 A_{ij}\partial_i \varphi
\partial_j \varphi dx +
\frac{\varphi_L^+-\varphi_L^-}{2L}\int_{\Omega_{-L,L}}\eta^2\left(\rho(|\nabla\varphi_1|^2)\nabla\varphi_1-\rho(|\nabla\varphi_2|^2)\nabla\varphi_2\right)\cdot
\vec{e}_ndx
$$
\begin{equation}\label{n7}
= -2\int_{\Omega_{-L-1,-L}}\eta (\varphi -\varphi_L^-)
A_{ij}\partial_i\eta
\partial_j\varphi dx-2\int_{\Omega_{L,L+1}}\eta (\varphi -\varphi_L^+)
A_{ij}\partial_i\eta
\partial_j\varphi dx.\ \ \ \ \ \ \
\end{equation}
The second integral on the left hand side of (\ref{n7}) vanishes. Indeed,
$$
\left.\begin{array}{l}
 \displaystyle{~~~~\int_{\Omega_{-L,L}}\eta^2(\rho(|\nabla\varphi_1|^2)\nabla\varphi_1-\rho(|\nabla\varphi_2|^2)\nabla\varphi_2)\cdot
\vec{e}_ndx}\\[3mm]
  \displaystyle{=\int^L_{-L}\eta^2(t)\int_{S_t}\left(\rho(|\nabla\varphi_1|^2)\nabla\varphi_1\cdot\vec{e}_n-\rho(|\nabla\varphi_2|^2)\nabla\varphi_2\cdot
\vec{e}_n\right) dx'dt=0,}
   \end{array}
 \right.
$$
since the two solutions possess the same mass flux $m_0$, here $S_t
= \Omega\cap\{x_n=t\}$ for $t\in[-L,L]$.

It follows from (\ref{n7}) and (\ref{n6}) that
\[\ba
&\ \ \ \ \lambda\int_{\Omega_{-L,L}}|\nabla\varphi|^2 dx\\
&\leq 4\Lambda\int_{\Omega_{-L-1,-L}}|\varphi-\varphi^-_L|
|\nabla\varphi| dx +
4\Lambda\int_{\Omega_{L,L+1}}|\varphi-\varphi^+_L| |\nabla\varphi|
dx\\
&\leq 2\Lambda\left(\int_{\Omega_{-L-1,-L}}|\varphi-\varphi^-_L|^2
dx + \int_{\Omega_{L,L+1}}|\varphi-\varphi^+_L|^2 dx+
\int_{\Omega_{-L-1,-L}\cup\ \Omega_{L,L+1}}|\nabla\varphi|^2
dx\right) \ea
\]
Due to the uniform Poincar\'{e} inequality, ie.
$$
\int_{\Omega_{-L-1,-L}}|\varphi-\varphi^-_L|^2 dx\leq C
\int_{\Omega_{-L-1,-L}}|\nabla\varphi|^2
dx,$$and$$\int_{\Omega_{L,L+1}}|\varphi-\varphi^+_L|^2 dx\leq C
\int_{\Omega_{L,L+1}}|\nabla\varphi|^2 dx,
$$
where $C$ is independent of $L$, we have
\begin{equation}\label{n8}
\lambda\int_{\Omega_{-L,L}}|\nabla\varphi|^2 dx \leq C
\int_{\Omega_{-L-1,-L}\cup\ \Omega_{L,L+1}}|\nabla\varphi|^2 dx.
\end{equation}
By the estimate (\ref{L}), one has
$$
\int_{\Omega_{-L-1,-L}\cup\ \Omega_{L,L+1}}|\nabla\varphi_k|^2 dx\leq C m_0^2\left(|\Omega_{-L-1,-L}|+|\Omega_{L,L+1}|\right),\ \ \ k=1,2,
$$
which implies that
$$
\int_{\Omega_{-L-1,-L}\cup\ \Omega_{L,L+1}}|\nabla\varphi|^2 dx\leq C,
$$
where $C$ is independent of $L$.

Combining this with (\ref{n8}) shows
$$
\int_{\Omega_{-L-1,-L}}|\nabla\varphi|^2 dx \rightarrow 0\ \ \ \text{and}\ \ \  \int_{\Omega_{L,L+1}}|\nabla\varphi|^2 dx \rightarrow 0,\ \ \
\text{as}\ \ L\rightarrow\infty.
$$
Taking $L\rightarrow\infty$ in (\ref{n8}) yields
$$
\nabla\varphi=0\ \ \ \ in \ \ \Omega.
$$
\endproof
As a direct application of the uniqueness, we can obtain the explicit form of the subsonic solution $\fai(x)$ to the Problem \ref{op}, provided that
the nozzle is a cylinder.

\begin{corollary}\label{lem-asym1} (Cylinder case)
Suppose that $\O$ is a cylinder, that is, $\O =
S\times(-\infty,+\infty)$, $S$ is a $n-1$ dimensional, simply
connected, $C^{2,\a}$ domain. Then the unique solution to Problem \ref{op} is given by
\[
\fai=q_0x_n+\fai_0,
\]
where $\fai_0$ is an arbitrary constant, $q_0$ is a constant defined
by
\[
\rho(q_0^2)q_0 = \f{m_0}{|S|}.
\]
\end{corollary}
\begin{figure}[!h]\centering
\includegraphics[width=110mm]{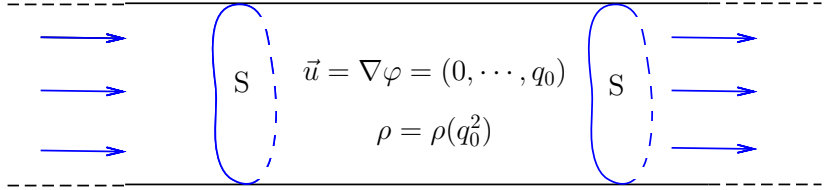}\\
\caption{Subsonic flow in cylinder case} \label{fig3}
\end{figure}

\section{Existence of the critical incoming mass flux}
In the Section 3 and Section 4, we have obtained the
existence of the uniformly subsonic flows associated with suitable
small incoming mass flux $m_0$ and the uniqueness of the uniformly
subsonic flow. In the following, it will be shown that there exists
a critical mass flux $M_c$ such that the flow is always uniformly
subsonic, provided that the mass flux $m_0$ is less than $M_c$.
\begin{theorem}\label{61}
Suppose the nozzle satisfies the basic assumptions (\ref{H1}).
Then there exists a positive constant $M_c \leq 1$, which depends
only on $\Omega$, such that if $0\leq m_0 < M_c$, then the following
problem
$$\left
 \{\begin{array}{ll}
 \displaystyle{\text{div} \left(\rho(|\nabla \varphi|^2)\nabla \varphi\right)=0},\ \ \ \ \ \ & \text{in}\ \ \ \Omega,\\
  \displaystyle{\frac{\partial\varphi}{\partial{\vec n}}=0,}& \text{on}\ \ \
  \p\O,\\
  \displaystyle{\int_S\rho\left(|\nabla\varphi|^2\right)\frac{\partial\varphi}{\partial \vec{l}}dS=m_0}
 \end{array}
 \right.
$$
has a unique uniformly subsonic solution $\fai(x)$ up to a constant
satisfying
$$\displaystyle{Q(m_0)=\sup_{x\in\bar\Omega}|\nabla\varphi|<1}.$$
Moreover, $Q(m_0)$ ranges over $[0, 1)$ as $m_0$ varies in $[0,
M_c)$.
\end{theorem}
\proof Choosing a strictly increasing sequence
$\displaystyle{\{q_n\}_{n=1}^\infty}$ satisfying
$\displaystyle{\lim_{n\rightarrow\infty}q_n=1}$. Consider the
following truncated problem
\begin{equation}\label{n9}
\left
 \{\begin{array}{ll}
 \displaystyle{\text{div} \left(\rho_n(|\nabla \varphi|^2)\nabla \varphi\right)=0},\ \ \ \ \ \ & \text{in}\ \ \ \Omega,\\
  \displaystyle{\frac{\partial\varphi}{\partial{\vec n}}=0,}& \text{on}\ \ \
  \p\O,\\
  \displaystyle{\int_S\rho_n\left(|\nabla\varphi|^2\right)\frac{\partial\varphi}{\partial
  \vec{l}} dS=m},
 \end{array}
 \right.
\end{equation}
where
$$\rho_n(s)=\left
 \{\begin{array}{ll}
 \displaystyle{\rho(s)},\ \ \ \ \ \ & \text{if}\ \ 0\leq s \leq q_n^2,\\[2mm]
 \text{smooth and decreasing}\ \ \ & \text{if}\ \ q_n^2\leq s\leq
 \left(\frac{q_n+1}2\right)^2,\\
  \rho\left(\left(\frac{q_n+1}{2}\right)^2\right),\ \ \ \ \ \ \ \ \ & \text{if}\ \ s > \left(\frac{q_n+1}2\right)^2,\\
   \end{array}
 \right.
$$
satisfies $\rho_n(s) +2s\rho_n'(s) > \lambda_n$ with some $\lambda_n
> 0$ for all $s \geq 0$. Let $\varphi_n(\cdot; m)$ solve the problem
(\ref{n9}) and set
$$\displaystyle{Q_n(m) =\sup_{x\in
\bar\Omega}|\nabla\varphi_n(\cdot; m)|}.$$

We claim that $Q_n(m)$ is a continuous function of $m$.

In fact, we take a sequence $m_j\rightarrow m$, it suffices to
prove
$$
\displaystyle{|\nabla\varphi_n(\cdot; m_j)|\rightarrow
|\nabla\varphi_n(\cdot; m)|}.
$$
Without loss of generality, we assume that there exists a positive constant $\bar M$, such that
$\displaystyle{\sup_{j\geq 1}m_j<\bar M}$.

It follows from Section 3 that the solution
$\varphi_n(\cdot; m_j)$ to the problem (\ref{n9}) with the mass flux
$m_j$ satisfies the H\"older gradient estimate
\begin{equation}\label{n10}
\parallel\nabla\varphi_n(\cdot;
m_j)\parallel_{C^{1,\alpha}(\Omega)}\leq C(\bar M)
\end{equation}
and
\begin{equation}
\parallel\varphi_n(\cdot;
m_j)-\varphi_n(0; m_j)\parallel_{C^{2,\alpha}(\Omega_L)}\leq C(\bar M,
L)\ \ \ \ \ \ for\ \ any\ \  L>0.
\end{equation}
Therefore, by Arzela-Ascoli Lemma and a diagonal argument, there
exists a subsequence $\varphi_n(\cdot; m_{j_k}) -\varphi_n(0;
m_{j_k})$ such that for any $L > 0$ and $0< \beta < \alpha$
$$
\left(\varphi_n\left(\cdot; m_{j_k}\right) -\varphi_n\left(0;
m_{j_k}\right)\right)\rightarrow \varphi_n(\cdot)\ \ \ \ \
\text{in}\ \ \ C^{2,\beta}(\Omega_L)\ \ \ \text{as}\ \
m_{j_k}\rightarrow m.
$$
And $\varphi_n(\cdot)$ solves the boundary value problem (\ref{n9})
and satisfies that
$$
\parallel\nabla\varphi_n\parallel_{C^{1,\alpha}(\Omega)}\leq C(\bar M).
$$
On the other hand, it follows from the previous sections that there
exists a $\varphi_n(\cdot; m)$ which solves (\ref{n9}). We can conclude that
$$
\nabla\varphi_n(\cdot) = \nabla\varphi_n(\cdot; m)
$$
by the uniqueness.

Hence, for any $L > 0$
$$
\nabla\varphi_n(\cdot;m_j)\rightarrow\nabla\varphi_n(\cdot; m)\ \ \
\ in \ \ C^{1,\beta}(\Omega_L),\ \ \ \ as\ \ m_j\rightarrow m,
$$
which proves the claim.

It follows from the claim that, there exists the largest $Q_n >0$
and the smallest $S_n > 0$ such that
$$
q_{n-1} < Q_n(m) < q_n, \ \ \ \ \ \ \text{for}\ \ \ \text{any}\ \
m\in (m_n, M_n).
$$
Moreover, clearly $M_{n+1} \geq M_n$. Set $\displaystyle{M_c =
\lim_{n\rightarrow\infty}M_n}$. It follows the definition of $M_n$
that $M_n \leq \rho(Q_n^2(M_n))Q_n(M_n) <1$, hence $M_c \leq 1$.

Then we can conclude that there exists a critical mass flux $M_c
\leq 1$, for any $m_0 < M_c$, there is $M_n$ such that $M_n > m_0$,
then $$Q(m_0) = Q_n(m_0) < q_n < 1.$$

Moreover, for any normalized subsonic speed $Q\in (0, 1)$, there
exists some $n$, such that $Q\in (0, q_n)$, therefore, there exists
a $m_0\in (0, M_n)$, such that $Q(m_0)=Q_n(m_0)=Q$ by the continuity
of $Q_n(m)$.

This completes the proof of Theorem \ref{61}.
\endproof

\section{Properties of the subsonic flow}
In this section, we consider the asymptotic behavior of the
uniformly subsonic flows at the far fields under the asymptotic assumption (\ref{AA}).
\begin{proposition}\label{pn1} Suppose that the nozzle satisfies the asymptotic assumption (\ref{AA}). Then the subsonic flow constructed before approaches
to uniform flows
at the far fields, ie.
$$
\nabla \varphi = (0, \cdots, q_\pm),\ \ \ \ \ \ \ as\ \ \
x_n\rightarrow \pm\infty,
$$
respectively, $q_\pm$ are constants uniquely determined by
$$
\rho(q_\pm^2)q_\pm =\frac{m_0}{|S_\pm|},
$$
respectively.
\end{proposition}
\proof Assume that $\varphi(x)$ is a classical solution of
$$\left
 \{\begin{array}{ll}
 \displaystyle{\text{div} (\rho(|\nabla \varphi|^2)\nabla \varphi)=0},\ \ \ \ \ \ & \text{in}\ \ \ \Omega,\\[3mm]
  \displaystyle{\frac{\partial\varphi}{\partial{\vec{n}}}=0,}& \text{on}\ \ \
  \p\O,\\
 \displaystyle{\int_S\rho\left(|\nabla\varphi|^2\right)\frac{\partial\varphi}{\partial
 \vec{l}}dS=m_0},
 \end{array}
 \right.
$$
satisfying
\begin{equation}\label{n1}
\parallel\nabla\varphi\parallel_{C^{1,\alpha}(\Omega)} \leq Cm_0.
\end{equation}
{\bf Step 1. A Special Case. } Suppose that $\Omega\cap\{x_n\geq
L_0\}= U_+\times[L_0, +\infty)$ for some $L_0$. Define a sequence of
functions as follows
$$\displaystyle{\varphi_k{(x',
x_n)=\varphi(x', x_n+k)}\chi_{\Omega_k}},$$ here $\Omega_k
=\{(x',x_n)|(x', x_n+k)\in \Omega,\ x_n+k> L_0+1\}$.

For any compact set $S\subset \overline{S}_+$ and $k$ sufficiently
large, it follows from the gradient estimate (\ref{n1}) that
$$
\parallel\nabla\varphi_k\parallel_{C^{1,\alpha}(S\times[-k/2, k/2])} \leq
C,
$$
where $C$ is independent of $k$. Set $\hat\varphi_k(x) =
\varphi_k(x) - \varphi_k(0)$,  for any fixed $L\geq 1$, if $k > 2L$,
we have
$$
\parallel\hat\varphi_k\parallel_{C^{2,\alpha}(S\times[-L, L])} \leq C,
$$
with $C$ independent of $k$. Therefore, by Ascoli-Arzela Lemma
and a diagonal procedure, there exists a subsequence
$\hat\varphi_{k_j}$, such that for any $L$
$$
\hat\varphi_{k_j}\rightarrow \varphi_0,\ \ \ \ \ \text{in} \ \
C^{2,\beta}(S\times[-L, L]) \ \ \ \text{with}\ \ \beta<\alpha,
$$
for any compact set $S\subset \overline{S}_+$.  Therefore
$\varphi_0$ solves the following problem
\begin{equation}\label{n2}
\left
 \{\begin{array}{ll}
  \displaystyle{\text{div} (\rho(|\nabla \varphi|^2)\nabla \varphi)=0},\ \ \ \ \ \ & \text{in}\ \ \ \ E_+=S_+\times(-\infty,+\infty),\\
  \displaystyle{\frac{\partial\varphi}{\partial\vec{n}}=0,}&\text{ in}\ \ \ \ \partial
  S_+\times(-\infty, +\infty),\\
  \displaystyle{\int_S\rho\left(|\nabla\varphi|^2\right)\frac{\partial\varphi}{\partial
  \vec{l}}dS=m_0},\\
  \varphi(x) = 0, & \text{on}\ \ \ x_n=0.
 \end{array}
 \right.
\end{equation}
Moreover,
$$
\nabla\varphi_k=\nabla\hat\varphi_k\rightarrow \nabla \varphi_0\ \ \
\text{in} \ \ C^{1,\mu}(S\times[-L, L])\ \ \ \text{for}\ \ \mu <
\beta.
$$
So, choosing $S=\overline{S}_+$ and $L=2$, we have
$$
\parallel\nabla \varphi_k-\nabla\varphi_0\parallel_{C^{\mu}(\overline{S}_+\times
[-2, 2])}\rightarrow 0\ \ \ \ \ \ \  \text{as}\ \ k\rightarrow
+\infty.
$$
By the definition of $\varphi_k$ and Corollary \ref{lem-asym1}, it
follows that
$$
\nabla \varphi\rightarrow \nabla\varphi_0=(0, \cdots, q_+)\ \ \
\text{as}\ \ x_n\rightarrow +\infty.
$$
This completes the proof of Proposition \ref{pn1} in this special case.

{\bf Step 2. General Case.} Suppose now that the nozzle satisfies (\ref{AA}). we can also define a sequence of functions as
$$\displaystyle{\varphi_k{(x',
x_n)=\varphi(x', x_n+k)}\chi_{\Omega_k}},
$$
here $\Omega_k =\{(x',x_n)|(x', x_n+k)\in \Omega,\ x_n+k> 1\}$. Then
similar to the Step 1, we can show that \be\label{44}
\nabla\varphi_{k_j}\rightarrow \nabla\varphi_0\ \ \ \text{in} \ \
C^{1,\beta}(S\times[-L, L]) \ee for any compact set $S\subset S_+$
and any fixed $L$, here $S$ may not reach the boundary $\partial
S_+$, and $\varphi_0$ is still the solution of boundary value
problem (\ref{n2}).

In particular, $\varphi_0$ satisfies the no-flow boundary condition
on the nozzle wall. Indeed, for any given point $(y', y_n)\in
\partial S_+\times(-\infty, +\infty)$, $\vec n = (\vec n_1, 0)$ is
the outer normal direction of the cylinder $S_+\times (-\infty,
+\infty)$ at $(y', y_n)$. For any $\delta
> 0$, there exists suitable large $K_0>0$, such that  $$(y'-\delta \vec{n}_1, y_n + k)\in S\times\{x_n=y_n+k\}\ \ \ \text{for}\ \ \ k > K_0,$$
where $S$ is a compact set of $\O\cap\{x_n=y_n+k\}$.

There exists a sequence of $n-1$ dimensional vectors $\{\vec z_k\}_{k=1}^\infty$, such that
$(y'-\delta \vec n_1+\vec z_k, y_n+k)\in \partial \Omega$, and
$\displaystyle{|\vec z_k|= dist ((y'-\delta\vec n_1, y_n+k),
\partial \Omega)}$. $\vec n_k$ is the out normal of the domain $\Omega$ at
$(y'-\delta \vec n_1+\vec z_k, y_n +k)$. Obviously,
$$
\displaystyle{\lim_{k\rightarrow+\infty}|\vec z_k|\rightarrow 0},\ \
\ \text{and}\ \ \  \displaystyle{\lim_{k\rightarrow+\infty}\vec
n_k\rightarrow \vec n},
$$
due to the assumption (\ref{AA}) on the nozzle at the far fields.

Therefore
$$
\left.\begin{array}{rl}
  \displaystyle{\ \ \ \nabla \varphi_0(y' - \delta \vec n_1, y_n)\cdot \vec{n}
  }&
  \displaystyle{=\lim_{k\rightarrow+\infty}\nabla\varphi(y'-\delta \vec n_1, y_n+k)\cdot
  \vec{n}}\\
  &\displaystyle{=\lim_{k\rightarrow+\infty}\left(\nabla\varphi(y'-\delta \vec n_1, y_n+k)- \nabla\varphi(y'-\delta \vec{n_1}+\vec z_k,
  y_n+k)\right)\cdot
  \vec{n}}\\
   &\displaystyle{\ \ \ + \lim_{k\rightarrow+\infty}\nabla\varphi(y'-\delta \vec n_1+ \vec z_k, y_n+k)\cdot
  (\vec{n}-\vec n_k)}\\
  \ &\displaystyle{=0}.
    \end{array}\right.
$$
As a consequence, $\displaystyle{\frac{\partial
\varphi_0}{\partial \vec n}(y', y_n)= 0}$.

Set
$$
\O'(\delta)=\left\{\left.x'\in \overline \O\cap\{x_n=k\}\right| dist(x',
\partial \O)<\delta\right\},\ \ \ \O_0'(\delta) = \left(\O\cap\{x_n=k\}\right)\backslash
\O'(\delta),
$$
and
$$
B_{i,\delta}(y_i')=\left\{\left.x'\in\mathbb{R}^{n-1}\right||x'-y_i'|<2\delta\right\},
\ \ \ y_i'\in \partial \O\cap\{x_n=k\},\ \ i=1, \cdots, N,
$$
such that
$$
\O'(\delta) \subset \bigcup_{i=1}^N B_{i,\delta}(y_i),\ \ \
\text{for\ \ any}\ \ \delta > 0.
$$
For any fixed $\delta > 0$, there exists a sufficiently large $K_0 >
0$, such that $\O_0'(\delta)\subset S$ for $k> K_0$, $S$ is a
compact set in $S_+$.

Choosing $L =2$, one has from (\ref{44}) that
\begin{equation}\label{n3}
\parallel \nabla \varphi -
\nabla\varphi_0\parallel_{C^{\mu}\left(\O_0'(\delta)\times[k,k+2]\right)}\rightarrow
0\ \ \ as\ \ k\rightarrow +\infty, \ \ \ \text{for} \ \ \mu<\beta.
\end{equation}
Near the boundary $\partial \O\cap \{x_n=k\}$, $\varphi$ possesses the following
estimates
$$
\parallel\nabla\varphi\parallel_{C^\eta \left(B_{i,\delta}^+(y_i')\times
(k,k+2)\right)}\leq C.
$$
for $\eta> 0$, with $\displaystyle{B_{i,\delta}^+(y_i')=B_{i,\delta}(y_i')\cap \left(\bar\O\cap\{x_n=k\}\right)}$.
Hence, for any $x', y'\in B_{i,\delta}^+(y_i')$, one
has
$$
|\nabla \varphi(x',k) - \nabla\varphi(y', k)|< C \delta^\eta.
$$
Then, for any $\varepsilon
> 0$, there exists $\delta >0$, such that
\begin{equation}\label{n4}
|\nabla \varphi(x',k) - \nabla\varphi(y', k)|<
\frac\varepsilon{N+1},\ \ \ \ \ \ \text{for \ \ any}\ \ x',\ y'\in
B_{i,\delta}^+(y_i'),
\end{equation}
and  $i= 1, 2, \cdots, N$.

On the other hand, it follows from (\ref{n3}) that there exists $K
> 0$ such that
\begin{equation}\label{n39}
|\nabla\varphi-\nabla\varphi_0|\leq \frac\varepsilon{N+1}, \ \ \ \ \
\ \ \text{for\ \ any}\ \ x\in \O_0'(\delta)\times(K, +\infty).
\end{equation}
Then, combining that (\ref{n4}) and (\ref{n39}), one can conclude
$$
|\nabla\varphi-(0,\cdots, q_+)| < \varepsilon,\ \ \ \ \ \ \text{for
\ \ any}\ \  x_n > K.
$$

Similarly, one can get the asymptotic behavior as
$x_n\rightarrow-\infty$. This completes the proof of
Proposition \ref{pn1}.
\endproof

{\bf Acknowledgements}
This existence theory in this paper is part of the PhD. thesis of Wei Yan written under the supervision of Zhouping Xin at the Chinese University of Hong Kong \cite{yan}. Parts of this work were done when Lili Du was a postdoctoral fellow in the Institute of Mathematical Science, the Chinese University of Hong Kong during Aug. 2008--Aug. 2009, he would like to thank the institute's support and hospitality.
Du is supported in part by NNSF of China (No. 10801055) and
SRFDP(No. 200805611026). Xin is supported in part by Zheng Ge Ru
Foundation, Hong Kong RGC Earmarked Research Grants CUHK 4040/06P, CUHK 4042/08P, and a Focus Area Grant from the Chinese University of Hong Kong. Yan is supported in part by the NNSF of China (No. 11071025),
the Foundation of CAEP (No. 2010A0202010) and Foundation of STCPL.

% BibTeX users please use
%\bibliographystyle{spmpsci}


\end{document}